\documentclass[reqno, a4paper, 11pt]{amsart}
\usepackage{amsmath,amsthm,amssymb}
\usepackage{graphicx}
\usepackage{time} 
\usepackage{comment} 
\usepackage{thmtools}
\usepackage{thm-restate} 

\usepackage[dvipdfmx]{hyperref}
\hypersetup{
	setpagesize=false,
	bookmarksnumbered=true,%
	bookmarksopen=true,%
	colorlinks=true,%
	linkcolor=magenta,
	citecolor=blue,
}
\theoremstyle{plain}
\newtheorem{theorem}{Theorem}[section]
\newtheorem{lemma}[theorem]{Lemma}
\newtheorem{corollary}[theorem]{Corollary}
\newtheorem{proposition}[theorem]{Proposition}

\theoremstyle{definition}
\newtheorem{definition}[theorem]{Definition}
\newtheorem{notation}[theorem]{Notation}
\newtheorem{remark}[theorem]{Remark}
\newtheorem{example}[theorem]{Example}
\newcommand{\CC}{\mathbb{C}}

\newcommand{\calD}{\mathcal{D}}
\newcommand{\Int}{\mathrm{int\,}}

\title[Splittable relative trisections of exotic 4-manifolds]{Splittable relative trisections of \\ exotic 4-manifolds}
\author[Natsuya Takahashi]{Natsuya Takahashi}
\date{November 28, 2026.}
\address{Department of Mathematics, College of Science and Technology, Nihon University, 1-8-14, Kanda-Surugadai, Chiyoda-ku, Tokyo 101-8308, Japan}
\email{takahashi.natsuya@nihon-u.ac.jp}
\begin{document}
\begin{abstract}
We consider the notion of a splittable trisection, whose associated diagram simply splits into two subdiagrams.
We construct splittable relative trisections for infinitely many exotic pairs of $4$-manifolds, with one piece representing a cork and the other a $4$-ball.
We relate exotic pairs obtained by cork twists via a cut-and-paste operation on trisection diagrams.
As a byproduct, we obtain relative trisections for an exotic pair whose induced open books on the boundaries are very similar.
\end{abstract}
\maketitle

\section{Introduction}

An effective method to construct exotic smooth structures on $4$-manifolds is to glue two $4$-manifolds with boundary in various ways.
Among such constructions, cork twisting, introduced by Akbulut~\cite{Akb91x}, provides a useful method. It consists of removing an embedded contractible $4$-manifold and regluing it in a different way.
In this paper, we study such constructions from the perspective of trisection theory.

A trisection, introduced by Gay--Kirby~\cite{GayKir16}, is a decomposition of a $4$-manifold into three pieces.
It gives a combinatorial description of a smooth structure on a 4-manifold via a trisection diagram, which consists of three systems of simple closed curves on a surface.
Trisections have been constructed for several exotic $4$-manifolds (see e.g. \cite{BaySae23}, \cite{ChaCotCui23x}, \cite{IslNay24}, \cite{Tak25b}, \cite{Tak25a}).

To describe a gluing of two $4$-manifolds with boundary, we consider the notion of a splittable trisection.
Roughly speaking, a trisection diagram is called splittable if there exists a collection of simple closed curves disjoint from all three curve systems such that cutting the diagram along this collection splits it into two relative trisection diagrams; see Definition~\ref{dfn:spl-rt}.
Although several gluings of relative trisections have been studied in \cite{CasGayPin18b}, \cite{CasOzb19}, \cite{KimMil20}, \cite{CasIslMilTom21}, \cite{GayMei22}, the gluing operation underlying our splitting is much simpler.
%
%
We show the existence of splittable trisections for exotic pairs.

\begin{restatable}{theorem}{mainexotictris}\label{main:exotictris}
There exist infinitely many exotic pairs $(Y_1,Y_2)$ of $4$-manifolds with boundary satisfying the following properties.
\begin{itemize}
\item Each $Y_i$ admits a splittable relative trisection diagram $\mathcal{D}_i$ with a single splitting circle.
\item Cutting $\mathcal{D}_i$ along the splitting circle yields two diagrams $\mathcal{D}'$ and $\mathcal{D}''$, which are independent of $i$.
\item The diagram $\mathcal{D}'$ represents a cork, while $\mathcal{D}''$ represents the $4$-ball.
\end{itemize}
\end{restatable}

\begin{figure}[!tbp]
\centering
\includegraphics[scale=0.25]{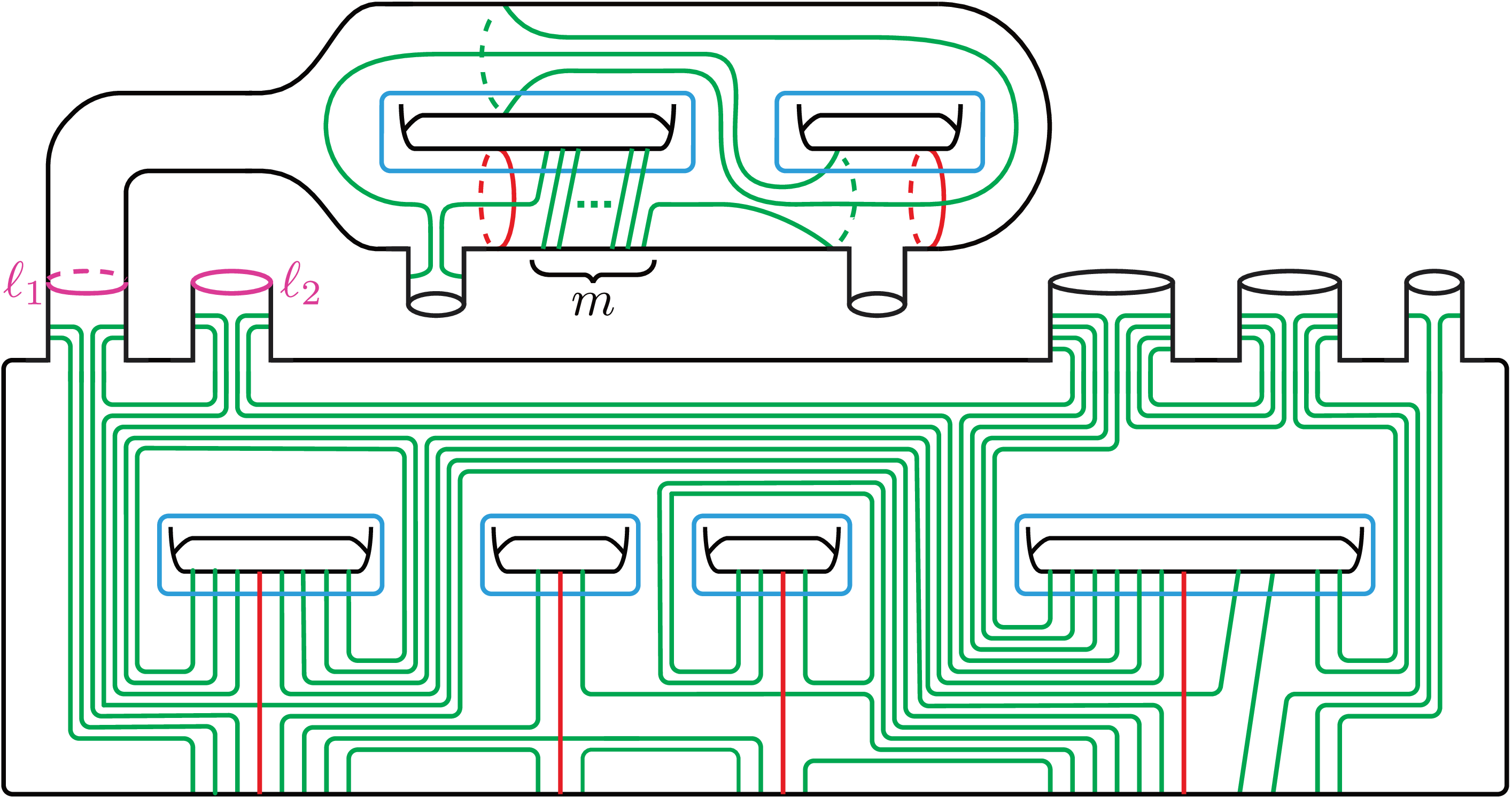}
\caption{A splittable trisection diagram of an exotic copy of the $4$-manifold shown in Figure~\ref{fig:td-exoticY2}.}
\label{fig:td-exoticY1}
\end{figure}

\begin{figure}[!tbp]
\centering
\includegraphics[scale=0.25]{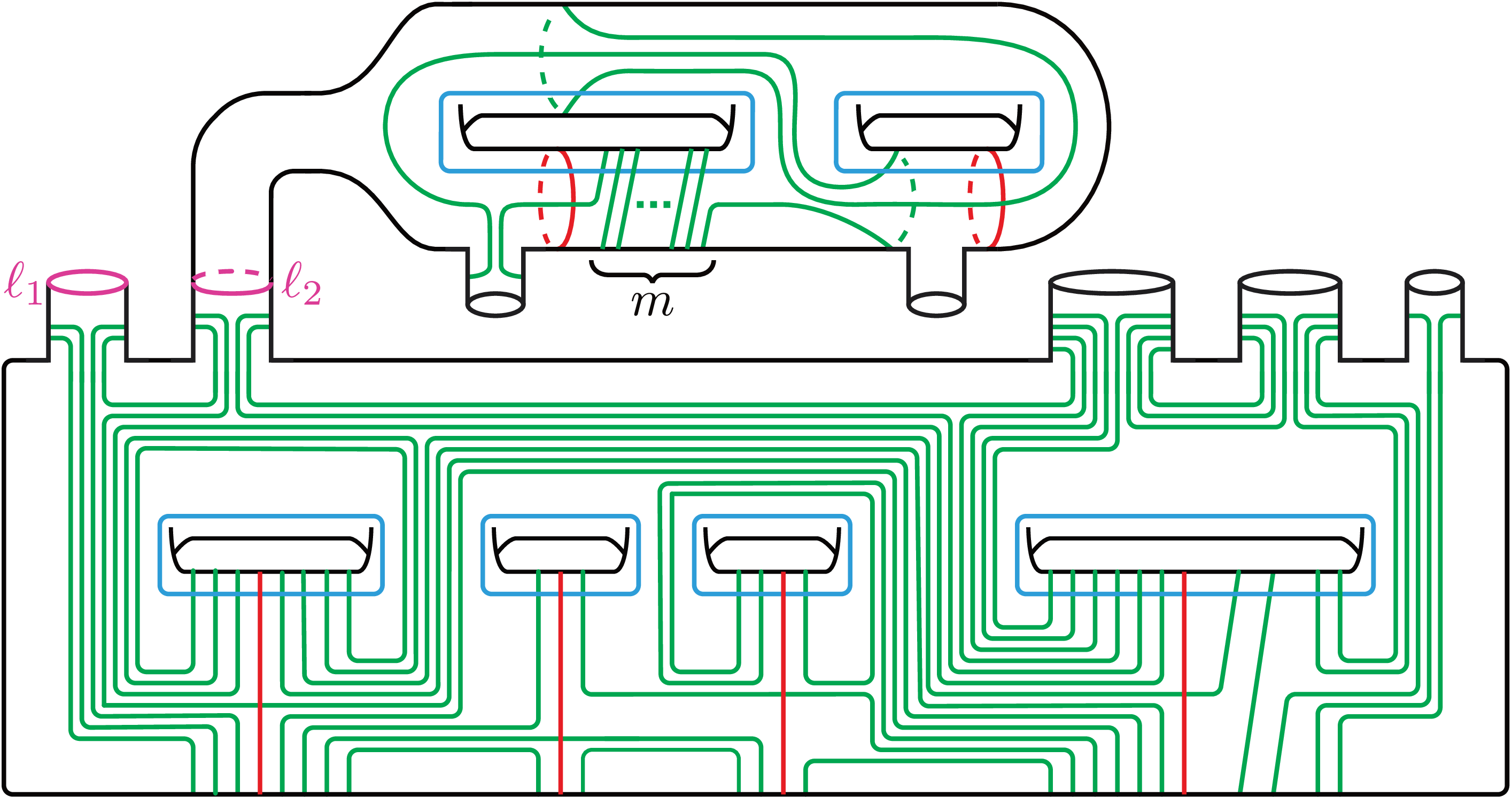}
\caption{A splittable trisection diagram of an exotic copy of the $4$-manifold shown in Figure~\ref{fig:td-exoticY1}.}
\label{fig:td-exoticY2}
\end{figure}

In other words, the differences in the smooth structures are encoded in how the two pieces of a relative trisection diagram are glued together.
Figures~\ref{fig:td-exoticY1} and~\ref{fig:td-exoticY2} show explicit examples of relative trisection diagrams for exotic pairs satisfying the conditions of Theorem~\ref{main:exotictris}.
In each diagram, the genus-$4$ piece below the splitting circle $\ell_i$ represents the Akbulut cork, while the genus-$2$ piece above $\ell_i$ represents $D^4$.
The latter piece corresponds to a $2$-handle attached to the Akbulut cork and thus represents a trisected $2$-handle.
The exotic pairs appearing in these examples are among those constructed by Akbulut--Yasui~\cite{AkbYas10}, \cite{AkbYas13}.
The following corollary is an immediate consequence of our examples.


\begin{corollary}
There exist infinitely many exotic pairs of $4$-manifolds with boundary such that both manifolds in each pair have trisection genus at most $6$.
\end{corollary}

Let $\mathcal{D}_A$ be the genus-$4$ subdiagram below $\ell_i$ in Figure~\ref{fig:td-exoticY1} (or~\ref{fig:td-exoticY2}), which represents the Akbulut cork and has distinguished boundary components $\ell_1$ and $\ell_2$.
We can construct pairs of homeomorphic $4$-manifolds by using the diagram $\mathcal{D}_A$.

\begin{restatable}{theorem}{maintdhomeo}\label{main:td-homeo}
Let $\mathcal{D}$ be a relative trisection diagram of a $2$-handlebody, where $\mathcal{D}$ has a distinguished boundary component $\ell$.
The two relative trisection diagrams obtained by gluing $\mathcal{D}_A$ and $\mathcal{D}$, identifying $\ell_i$ with $\ell$ for $i\in\{1,2\}$, give homeomorphic $4$-manifolds.
\end{restatable}

Furthermore, we study the open book decompositions on the boundary $3$-manifolds induced by splittable relative trisections.
We show that a relative trisection split along a single circle induces an open book that also splits into two pieces.
See Lemma~\ref{lem:srtd-obd} for details.
For each relative trisection of the exotic pairs shown in Figures~\ref{fig:td-exoticY1} and~\ref{fig:td-exoticY2}, the induced open book is obtained by gluing the monodromies induced by the two pieces along $\ell_i$: the genus-$4$ piece $\calD_A$ below $\ell_i$ and the genus-$2$ piece above $\ell_i$.

We also obtain relative trisection diagrams representing an exotic pair of $4$-manifolds whose induced open books are very similar.

\begin{corollary}\label{cor:obd-exotic}
There exist two relative trisection diagrams $\mathcal{D}_1$ and $\mathcal{D}_2$ with the same parameters $(g,k;p,b)$ that satisfy the following properties:
\begin{itemize}
\item The diagrams $\mathcal{D}_1$ and $\mathcal{D}_2$ represent an exotic pair of simply-connected $4$-manifolds whose boundaries are homology $3$-spheres.
\item The induced open books have the same page $P$, and their monodromies $\phi_1$ and $\phi_2$ satisfy $\phi_1\circ\phi_2^{-1}=\tau_{\ell_1}\circ\tau_{\ell_2}^{-1}$,
where $\ell_1$ and $\ell_2$ are simple closed curves each parallel to a boundary component of $P$, and $\tau_{\ell_i}$ denotes the right-handed Dehn twist along $\ell_i$.
\end{itemize}
\end{corollary}

In other words, there exists a diffeomorphism $f:P\to P$ such that $\phi_i=\tau_{\ell_i}\circ f$ for each $i$.
We give an explicit example of a pair of relative trisection diagrams satisfying these conditions, in which the diffeomorphism $f$ is induced by the genus-$4$ piece $\calD_A$ representing the Akbulut cork. 

\subsection*{Notation and conventions}

Throughout this paper, all manifolds and surfaces are assumed to be compact, connected, oriented, and smooth.
%
When we say that two manifolds are diffeomorphic, it implies that they are orientation-preservingly diffeomorphic.
%
%
%
Let $\Sigma_{g,b}$ denote a surface of genus $g$ with $b$ boundary components.

\section{Preliminaries}\label{sect:prel}
\subsection{Corks and exotic 4-manifolds}\label{subsec:corks}

We give an overview of corks and introduce exotic pairs obtained by cork twists, which will be used to construct relative trisections.
A pair $(C,\tau)$ is called a \textit{cork} if $C$ is a compact, contractible, smooth $4$-manifold and $\tau:\partial C\to\partial C$ is a diffeomorphism that extends to a self-homeomorphism of $C$, but cannot extend to any self-diffeomorphism of $C$.

\begin{figure}[!tbp]
\centering
\includegraphics[scale=0.23]{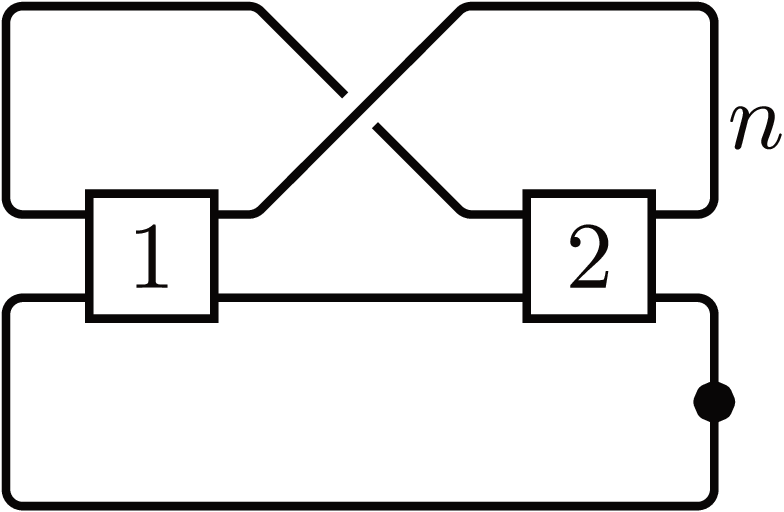}
\caption{A handlebody diagram of $C_n$.}
\label{fig:kd-Cn}
\end{figure}

We now give an example of a cork.
Let $C_n$ be the contractible $4$-manifold given by the handlebody diagram shown in Figure~\ref{fig:kd-Cn}.
The manifold $C_0$ admits a cork structure, and is the first example of a cork constructed by Akbulut~\cite{Akb91x}.

\begin{example}
%
Let $\tau:\partial C_0\to\partial C_0$ be the smooth involution obtained by first surgering the embedded $S^1\times D^3\subset C_0$ to $D^2\times S^2$, and then surgering the other embedded $D^2\times S^2\subset C_0$ to $S^1\times D^3$.
Equivalently, $\tau$ is obtained by interchanging the ``dot'' and the ``$0$'' in the handlebody diagram.
Then $(C_0,\tau)$ is a cork.
\end{example}

Let $(C,\tau)$ be a cork, and let $X$ be a $4$-manifold containing $C$ as a smooth submanifold.
Define
\[ X' := (X-\mathrm{int}{C})\cup_\tau C. \]
Thus, $X'$ is obtained from $X$ by removing the interior of $C$ and regluing $C$ by $\tau$.
This operation is called a \textit{cork twist}.
Since $\tau$ extends to a homeomorphism of $C$, the manifold $X'$ is homeomorphic to $X$, but $X'$ need not be diffeomorphic to $X$.
Matveyev \cite{Mat96} and Curtis--Freedman--Hsiang--Stong \cite{CurFreHsiSto96} showed that any two simply-connected, closed, exotic $4$-manifolds are related by a cork twist.

We introduce pairs of homeomorphic $4$-manifolds with boundary obtained by twisting along the Akbulut cork $C_0$.
For this purpose, we fix some notation.

\begin{figure}[!tbp]
\centering
\begin{minipage}[b]{0.49\columnwidth}
	\centering
	\includegraphics[scale=0.23]{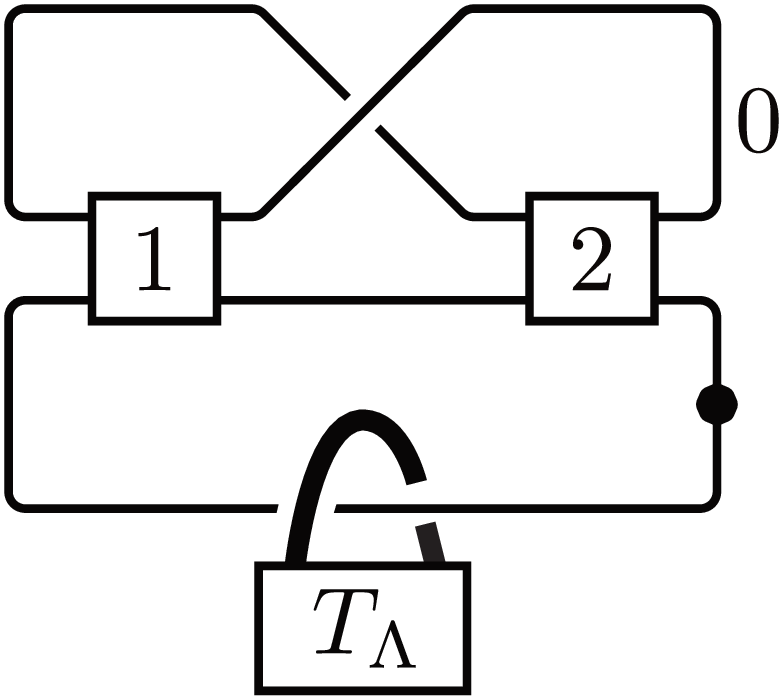}
	\caption{$Y_1(\Lambda)$}
	\label{fig:kd-Y1L}
\end{minipage}
\begin{minipage}[b]{0.49\columnwidth}
    \centering
	\includegraphics[scale=0.23]{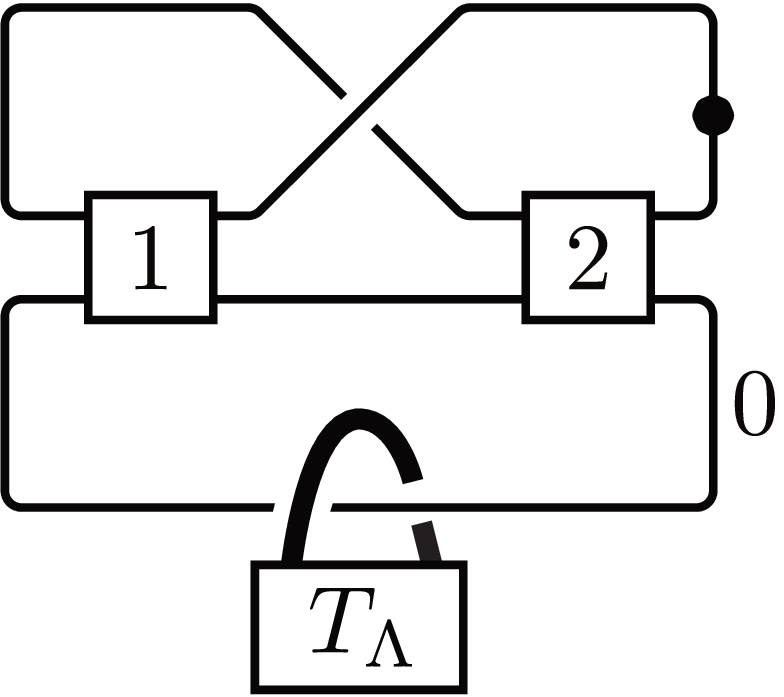}
	\caption{$Y_2(\Lambda)$}
	\label{fig:kd-Y2L}
\end{minipage}
\end{figure}

\begin{figure}[!tbp]
\centering
\includegraphics[scale=0.23]{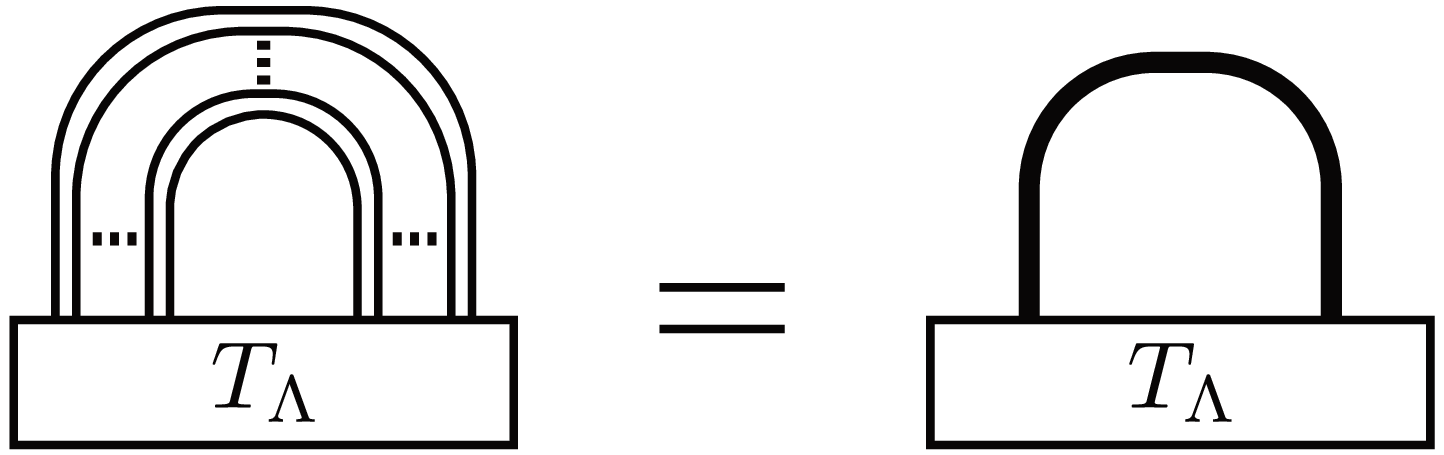}
\caption{$T_\Lambda$}
\label{fig:kd-tangle}
\end{figure}

\begin{notation}
Let $\Lambda=\{U_1,\ldots,U_s;K_1,\ldots,K_t\}$ denote a handlebody diagram of a $2$-handlebody, where each $U_i$ is a dotted circle and each $K_i$ is a framed knot.
Let $Y_1(\Lambda)$ and $Y_2(\Lambda)$ be the $4$-manifolds given by the handlebody diagrams in Figures~\ref{fig:kd-Y1L} and~\ref{fig:kd-Y2L}, respectively.
Here, the symbol $T_\Lambda$ in the box represents the tangle obtained from the diagram $\Lambda$, and the thick strands represent several parallel strands; see Figure~\ref{fig:kd-tangle}.
In particular, when $\Lambda=\{K\}$ represents the knot $m$-trace $X_m(K)$ and a single strand exiting the box passes through the cork, we write $Y_i(\Lambda)$ as $Y_i(K,m)$.
\end{notation}

For any $2$-handlebody $\Lambda$, the $4$-manifolds $Y_1(\Lambda)$ and $Y_2(\Lambda)$ are homeomorphic since they are related by twisting along $C_0$.
By the results of Akbulut--Yasui~\cite[Subsection~9.1]{AkbYas13} and~\cite[Theorem~3.3]{AkbYas10}, we obtain the following exotic pairs of $4$-manifolds with boundary.
Note that each of these papers gives a broader class of exotic pairs.
%

\begin{theorem}\label{thm:exoticAY}
Let $T_{p,p-1}$ be the $(p,p-1)$-torus knot.
Then $Y_1(T_{p,p-1},m)$ and $Y_2(T_{p,p-1},m)$ are exotic if $p=2$ and $m\leq -1$, or if $p\geq 3$ and $m\leq p^2-3p+1$.
\end{theorem}

\begin{figure}[!tbp]
\centering
\begin{minipage}[b]{0.49\columnwidth}
	\centering
	\includegraphics[scale=0.23]{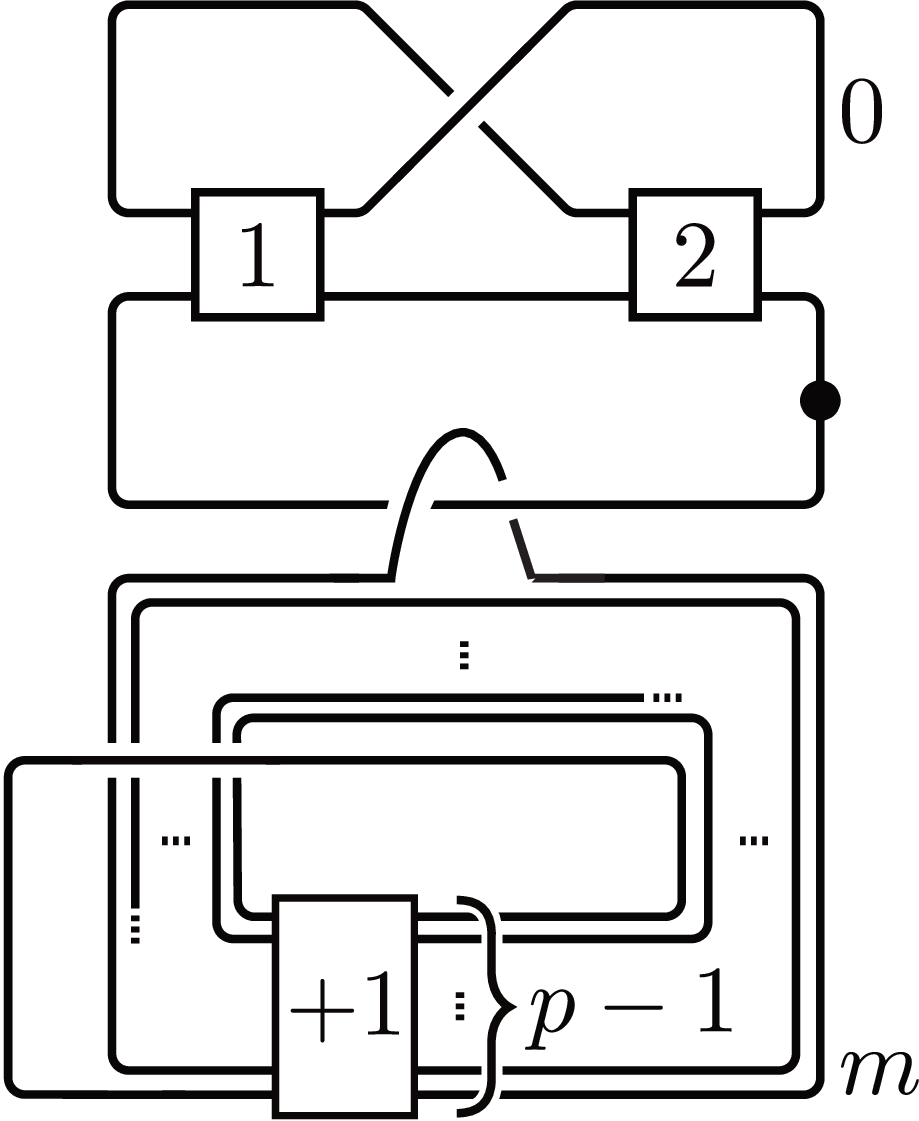}
	\caption{$Y_1(T_{p,p-1},m)$}
	\label{fig:kd-Ac+torus1}
\end{minipage}
\begin{minipage}[b]{0.49\columnwidth}
    \centering
	\includegraphics[scale=0.23]{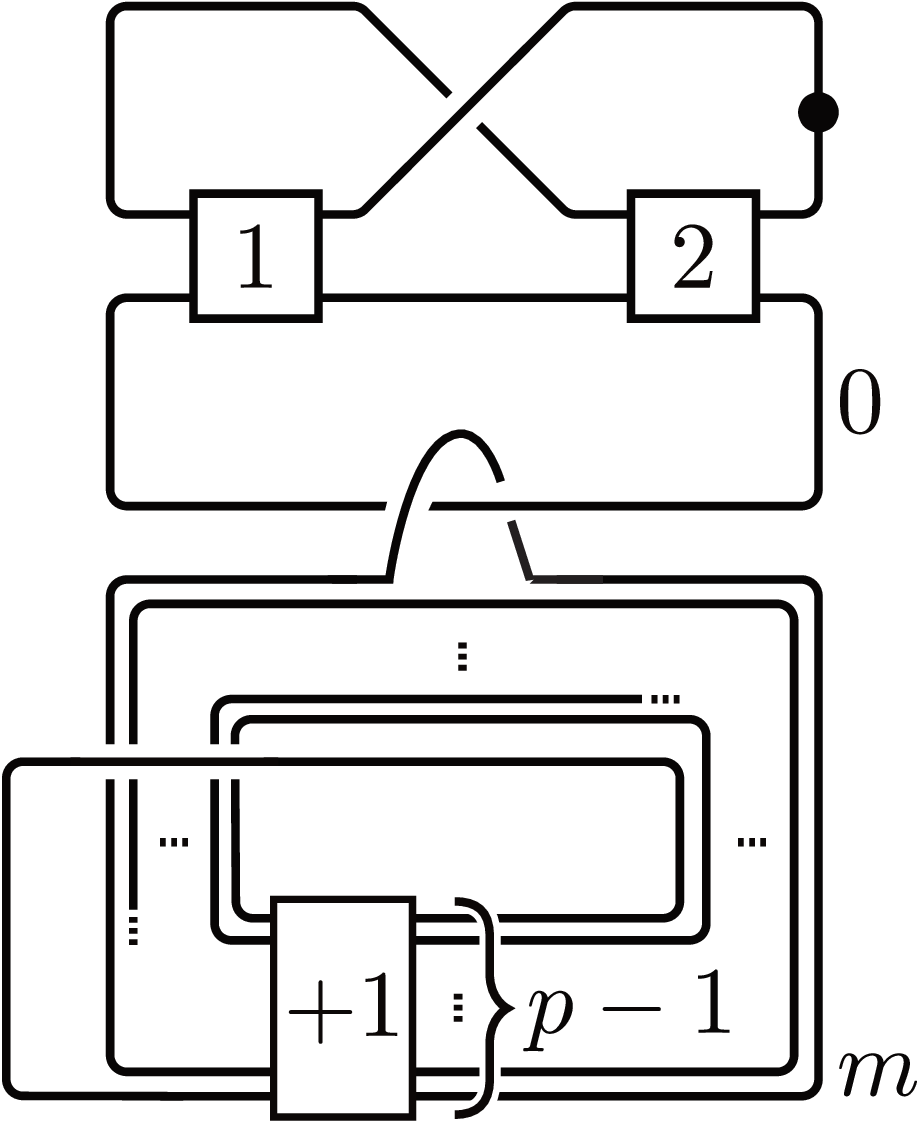}
	\caption{$Y_2(T_{p,p-1},m)$}
	\label{fig:kd-Ac+torus2}
\end{minipage}
\end{figure}

Note that $T_{p,p-1}$ is the unknot when $p=2$.
Theorem~\ref{thm:exoticAY} means that the $4$-manifolds given by Figures~\ref{fig:kd-Ac+torus1} and~\ref{fig:kd-Ac+torus2} form an exotic pair.
In Section~\ref{sect:tris-exotic}, we explicitly construct relative trisection diagrams for these exotic pairs.

\subsection{Relative trisections}\label{sec:tris}

We recall the definition of a relative trisection given by Castro--Gay--Pinz\'{o}n-Caiced~\cite{CasGayPin18a}.
The definition uses a model $Z_k$ of a $4$-dimensional $1$-handlebody equipped with a decomposition of its boundary, which we describe below.

\begin{figure}[!tbp]
\centering
\includegraphics[scale=0.23]{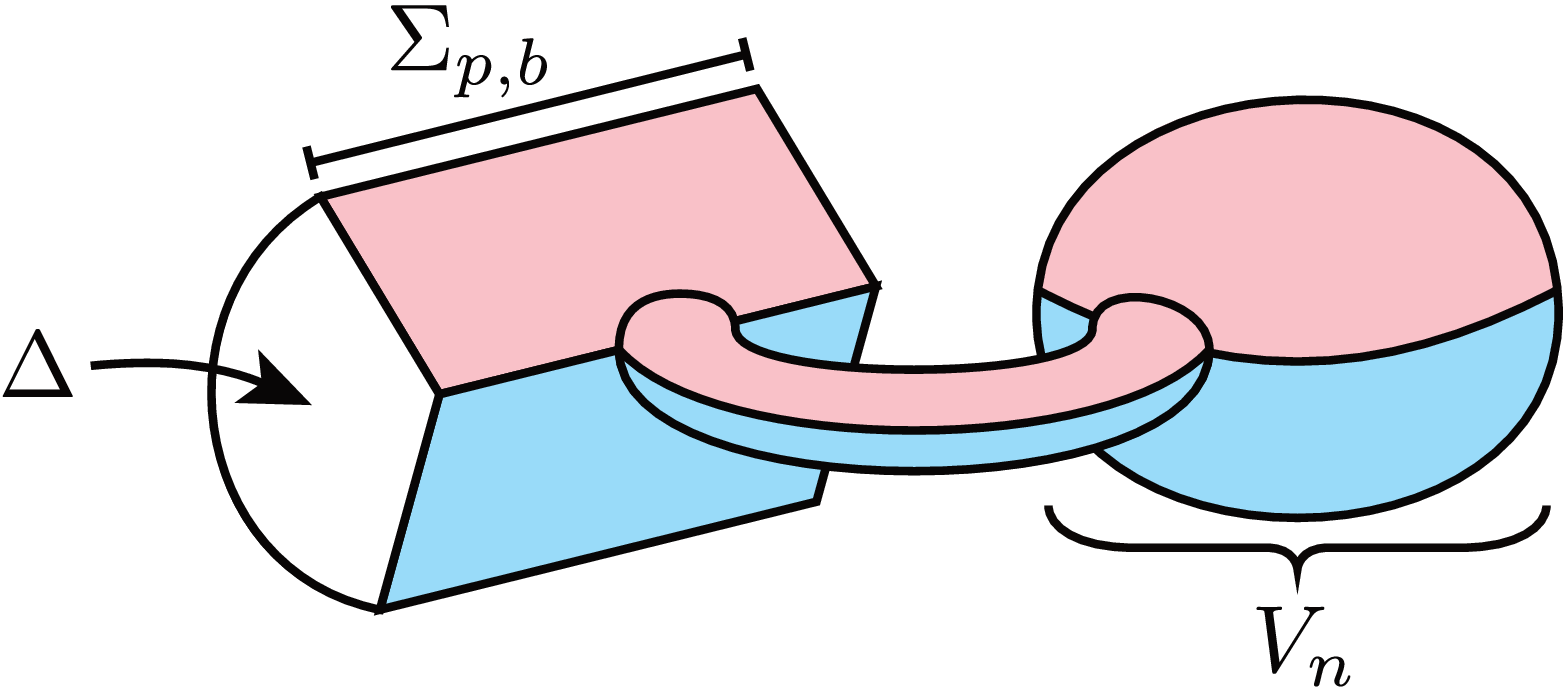}
\caption{A schematic picture of the model $Z_k$. The red and blue portions represent $Y^+_{g,k;p,b}$ and $Y^-_{g,k;p,b}$, respectively.}
\label{fig:modelZk}
\end{figure}

\begin{definition}
Let $\Delta$ be the sector of the unit disk defined by
\begin{equation*}
\Delta := \{z\in\CC \mid 0\leq |z|\leq 1,\ -\pi/3\leq \arg z\leq \pi/3\}.
\end{equation*}
We decompose its boundary as
$\partial\Delta=\partial^-\Delta\cup\partial^0\Delta\cup\partial^+\Delta$,
where $\partial^\pm\Delta$ are the radial edges corresponding to $\arg z=\pm\pi/3$, respectively, and $\partial^0\Delta$ is the remaining arc.
For integers $p\geq0$ and $b\geq1$, set $U_{p,b}:=\Sigma_{p,b}\times\Delta$ and decompose its boundary as
$\partial U_{p,b}=\partial^-U_{p,b}\cup\partial^0U_{p,b}\cup\partial^+U_{p,b}$,
where
\begin{align*}
\partial^\pm U_{p,b}:=\Sigma_{p,b}\times\partial^\pm\Delta
\quad\textrm{and}\quad
\partial^0 U_{p,b}:=(\Sigma_{p,b}\times\partial^0\Delta)\cup(\partial\Sigma_{p,b}\times\Delta).
\end{align*}
Let $g$ and $k$ be non-negative integers with $n:=k-2p-b+1\geq 0$.
Set $V_n:=\natural^n(S^1\times D^3)$, and equip its boundary with a genus-$(g-p)$ Heegaard splitting $\partial V_n=\partial^-V_n\cup\partial^+V_n$.
%
Note that $\#^n(S^1 \times S^2)$ has a unique Heegaard splitting up to isotopy in each genus.
We define $Z_k:=U_{p,b}\natural V_n$, where the boundary connected sum is taken along small $3$-balls centered at points in $\mathrm{int}(\partial^-U_{p,b}\cap\partial^+U_{p,b})$
and $\partial^-V_n\cap\partial^+V_n$.
See Figure~\ref{fig:modelZk}.
The boundary of $Z_k$ inherits a decomposition
\begin{equation*}
\partial Z_k=Y^-_{g,k;p,b}\cup Y^0_{p,b}\cup Y^+_{g,k;p,b},
\end{equation*}
where
\begin{align*}
Y^\pm_{g,k;p,b}:=\partial^\pm U_{p,b}\natural\partial^\pm V_n
\quad\textrm{and}\quad
Y^0_{p,b}:=\partial^0 U_{p,b}.
\end{align*}
\end{definition}

\begin{definition} 
Let $g$, $k$, $p$, and $b$ be integers with $g,k,p\geq0$, $b\geq1$, and $2p+b-1\leq k\leq g+p+b-1$.
Let $X$ be a compact, connected, oriented, smooth $4$-manifold with connected boundary.
A decomposition $X=X_1\cup X_2\cup X_3$ is called a $(g,k;p,b)$-\textit{relative trisection} if 
for each $i\in\{1,2,3\}$, there exists a diffeomorphism $\varphi_i:X_i\to Z_k$ such that
\begin{equation*}
\varphi_i(X_i\cap X_{i\pm1})=Y^\mp_{g,k;p,b}
\quad\textrm{and}\quad
\varphi_{i}(X_i\cap \partial{X})=Y^0_{p,b}.
\end{equation*}
where the indices are taken modulo $3$.
\end{definition}

By definition, the following properties hold.

\begin{itemize}
\item Each sector $X_i$ is diffeomorphic to $\natural^{k}(S^1\times D^3)$.
\item For distinct $i,j\in\{1,2,3\}$, the double intersection $X_i\cap X_j(=\partial X_i\cap\partial X_j)$ is diffeomorphic to a compression body from $\Sigma_{g,b}$ to $\Sigma_{p,b}$.
%
\item The triple intersection $X_1\cap X_2\cap X_3$ is diffeomorphic to $\Sigma_{g,b}$.
\item For each $i\in\{1,2,3\}$, $X_i \cap\partial{X}$ is diffeomorphic to $\Sigma_{p,b}\times [0,1]$.
\end{itemize}

By the fourth property, a $(g, k; p, b)$-relative trisection of a 4-manifold $X$ induces an open book decomposition on the boundary $\partial{X}$ with page $\Sigma_{p,b}$. 
For the induced open book, the binding is $\partial(X_1\cap X_2\cap X_3)$ and the page is diffeomorphic to $X_i\cap X_j \cap\partial{X}$.
%
For any $4$-manifold $X$ with boundary and any open book decomposition of $\partial X$, there exists a relative trisection of $X$ that induces the given open book decomposition (see \cite{GayKir16}, \cite{Cas16}, \cite{CasGayPin18a}, \cite{Sae26x}).

A relative trisection is encoded by a relative trisection diagram, which is a $4$-tuple $(\Sigma;\alpha, \beta, \gamma)$ consisting of a surface and three collections of curves.
We recall the definition given by Castro--Gay--Pinzón-Caicedo~\cite{CasGayPin18a}.
We recall the definition of a relative trisection diagram due to Castro--Gay--Pinzón-Caicedo~\cite{CasGayPin18a}.
We first define several equivalence relations on tuples consisting of a surface and collections of curves.

\begin{definition}
Let $\Sigma$ and $\Sigma'$ be compact, connected, oriented surfaces with boundary.
For each $i\in\{1,\ldots,n\}$, let $\eta^i=\{\eta^i_1,\ldots,\eta^i_k\}$ and $\zeta^i=\{\zeta^i_1,\ldots,\zeta^i_k\}$ be collections of $k$ pairwise disjoint simple closed curves on $\Sigma$ and $\Sigma'$, respectively.
Consider the two $(n+1)$-tuples $\mathcal{D}:=(\Sigma;\eta^1,\ldots,\eta^n)$ and $\mathcal{D}':=(\Sigma';\zeta^1,\ldots,\zeta^n)$.
\begin{itemize}
\item Suppose that $\Sigma=\Sigma'$. We say that $\mathcal{D}$ and $\mathcal{D}'$ are \textit{isotopic} if, for each $i\in\{1,\ldots,n\}$, there exists an ambient isotopy $\{\varphi_t^i:\Sigma\to\Sigma\}_{t\in[0,1]}$ satisfying $\varphi_1^i(\eta^i)=\zeta^i$.
\item Suppose that $\Sigma=\Sigma'$. We say that $\mathcal{D}$ and $\mathcal{D}'$ are \textit{slide-equivalent} if, for each $i\in\{1,\ldots,n\}$, $\eta^i$ can be modified into a collection isotopic to $\zeta^i$ by a sequence of slides among the curves in $\eta^i$.
\item We say that $\mathcal{D}$ and $\mathcal{D}'$ are \textit{slide-diffeomorphic} if there exists an orientation-preserving diffeomorphism $f:\Sigma\to\Sigma'$ such that the diagram $(f(\Sigma);f(\eta^1),\ldots,f(\eta^n))$ is slide-equivalent to $\mathcal{D}'$.
Here, each $f(\eta^i)$ denotes the collection of curves $\{f(\eta^i_1),\ldots,f(\eta^i_k)\}$.
\end{itemize}
\end{definition}

\begin{figure}[!tbp]
\centering
\includegraphics[scale=0.23]{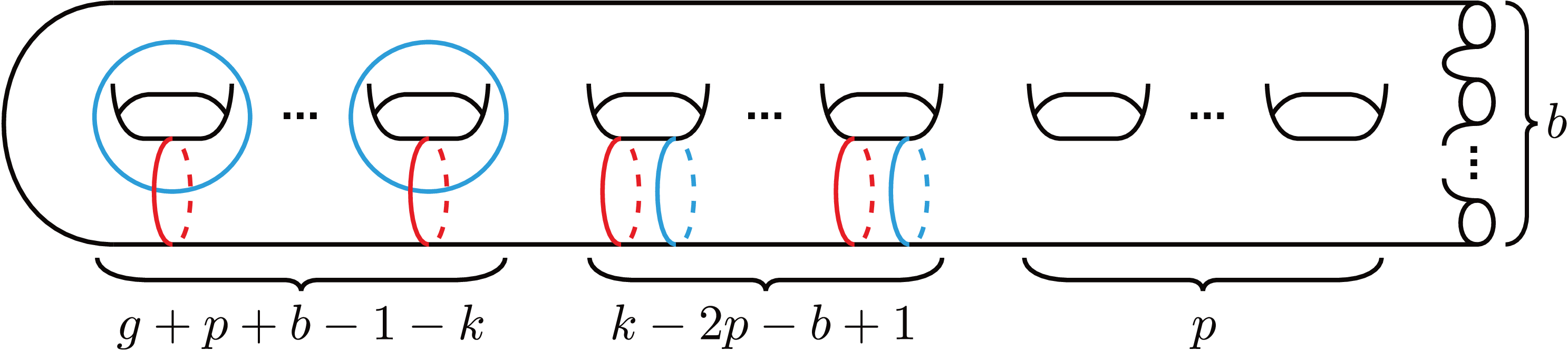}
\caption{The standard diagram $(\Sigma_{g,b}; \delta, \epsilon)$.}
\label{fig:td-gkpb}
\end{figure}

\begin{figure}[!tbp]
\centering
\includegraphics[scale=0.23]{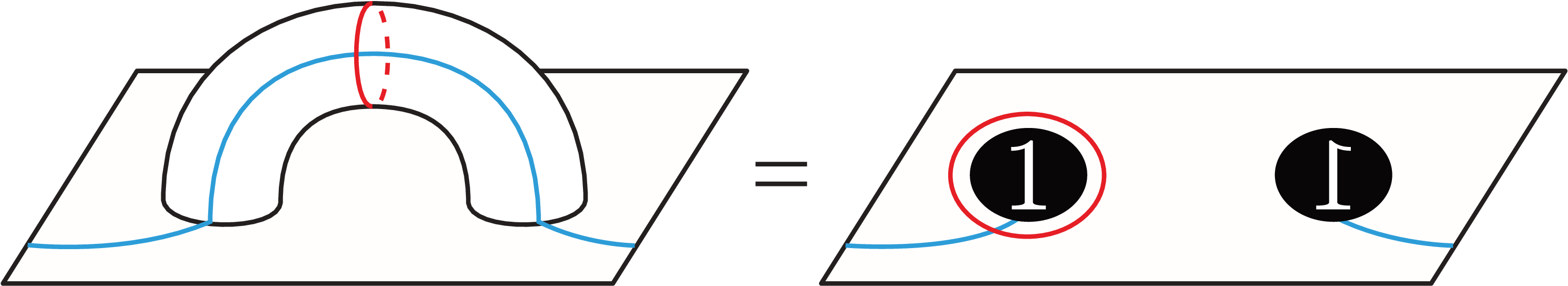}
\caption{A representation of a cylinder by a pair of labeled disks.}
\label{fig:diskpair}
\end{figure}

\begin{figure}[!tbp]
\centering
\includegraphics[scale=0.23]{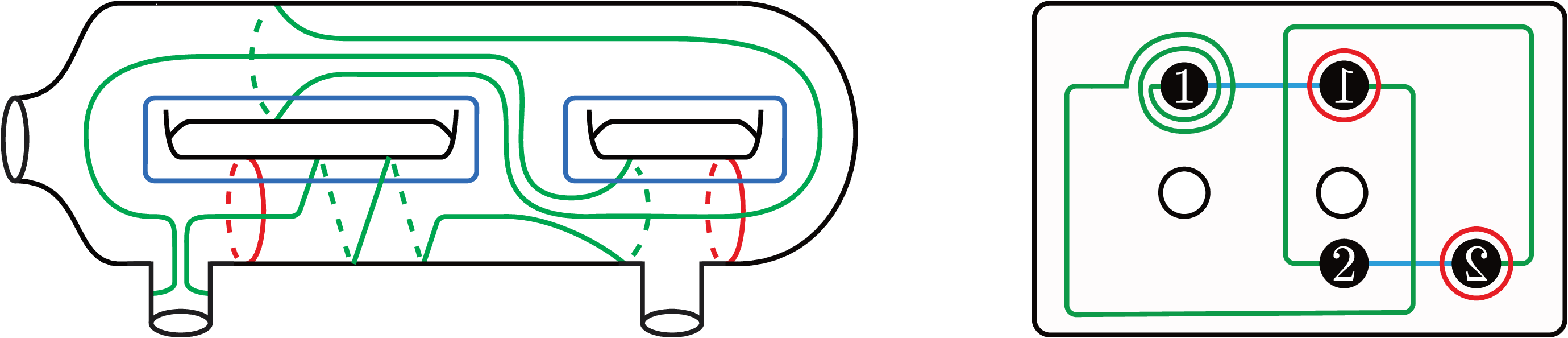}
\caption{Left: a $(3,2;0,3)$-relative trisection diagram. Right: a planar representation of the left diagram.}
\label{fig:td-exm}
\end{figure}

\begin{definition}\label{def}
Let $g$, $k$, $p$, and $b$ be integers satisfying $g,k,p\geq0$, $b\geq1$, and $2p+b-1\leq k\leq g+p+b-1$.
Let $\Sigma$ be a surface diffeomorphic to $\Sigma_{g,b}$, and let $\alpha=\{\alpha_1,\ldots,\alpha_{g-p}\}$, $\beta=\{\beta_1,\ldots,\beta_{g-p}\}$, and $\gamma=\{\gamma_1,\ldots,\gamma_{g-p}\}$ be collections of $g-p$ pairwise disjoint simple closed curves on $\Sigma$.
The $4$-tuple $\mathcal{D}=(\Sigma;\alpha,\beta,\gamma)$ is called a $(g,k;p,b)$-\textit{relative trisection diagram} if each of the three pairs $(\Sigma;\alpha,\beta)$, $(\Sigma;\beta,\gamma)$, and $(\Sigma;\gamma,\alpha)$ is slide-diffeomorphic to the standard diagram $(\Sigma_{g,b};\delta,\epsilon)$ depicted in Figure \ref{fig:td-gkpb}.
The collections of curves $\delta=\{\delta_1,\ldots,\delta_{g-p}\}$ and $\epsilon=\{\epsilon_1,\ldots,\epsilon_{g-p}\}$ are represented by the red and blue curves, respectively.
\end{definition}
For a relative trisection diagram $\mathcal{D}=(\Sigma;\alpha,\beta,\gamma)$, we refer to the curves in $\alpha$, $\beta$, and $\gamma$ as red, blue, and green, respectively.
We often depict relative trisection diagrams in the plane using the identification shown in Figure~\ref{fig:diskpair}. See, for example, Figure~\ref{fig:td-exm}.
%
%

There is a one-to-one correspondence between relative trisections and relative trisection diagrams.

\begin{theorem}[{\cite[Theorem~3]{CasGayPin18a}}]\label{thm:rt-rtd}
The following hold.
\begin{itemize}
\item For any $(g,k;p,b)$-relative trisection diagram $(\Sigma;\alpha,\beta,\gamma)$, there exists a unique (up to diffeomorphism) trisected $4$-manifold $X=X_1\cup X_2\cup X_3$ such that $X_1\cap X_2\cap X_3\cong\Sigma$ and, under this identification, the $\alpha$, $\beta$, and $\gamma$ curves bound compressing disks in $X_3\cap X_1$, $X_1\cap X_2$, and $X_2\cap X_3$, respectively.
\item Every relative trisection $X=X_1\cup X_2\cup X_3$ is induced by a relative trisection diagram under the above correspondence.
\item If two relative trisection diagrams $\mathcal{D}$ and $\mathcal{D}'$ induce diffeomorphic relative trisections, then $\mathcal{D}$ and $\mathcal{D}'$ are slide-diffeomorphic.
\end{itemize}
\end{theorem}

%

Recall that a $(g,k;p,b)$-relative trisection induces an open book decomposition on the boundary $3$-manifold with page $\Sigma_{p,b}$. 
There is an algorithm to compute the monodromy of the induced open book from a relative trisection diagram.

\begin{theorem}[\cite{CasGayPin18a}, \cite{CasOzb19}]\label{thm:rtd-obd}
Let $\calD = (\Sigma; \alpha, \beta, \gamma)$ be a $(g,k;p,b)$-relative trisection diagram. 
Let $\Sigma_{\alpha}$ be the surface obtained from $\Sigma$ by surgery along $\alpha$, and let $\mu_{\alpha} : (\Sigma-\alpha) \to \Sigma_{\alpha}$ be the induced map. 
Then we can take ordered collections $\mathfrak{a}$, $\mathfrak{b}$, $\mathfrak{c}$, $\mathfrak{a}'$, and $\mathfrak{a}^*$ of $2p + b - 1$ properly embedded arcs in the surface $\Sigma$ as follows:
\begin{enumerate}
\item Take a collection of properly embedded arcs $\mathfrak{a}$ in $\Sigma$ such that cutting $\Sigma_{\alpha}$ along $\mathfrak{a}$ yields a disk.
\item Take a collection of properly embedded arcs $\mathfrak{b}$ in $\Sigma$ and a collection of simple closed curves $\beta'$ such that $(\alpha,\mathfrak{b})$ is slide equivalent to $(\alpha,\mathfrak{a})$, $\beta'$ is slide equivalent to $\beta$, 
and $\mathfrak{b} \cap \beta' = \emptyset$.
\item Take a collection of properly embedded arcs $\mathfrak{c}$ in $\Sigma$ and a collection of simple closed curves $\gamma'$ such that $(\beta',\mathfrak{c})$ is slide equivalent to $(\beta',\mathfrak{b})$, $\gamma'$ is slide equivalent to $\gamma$, and $\mathfrak{c} \cap \gamma' = \emptyset$.
\item Take a collection of properly embedded arcs $\mathfrak{a}'$ in $\Sigma$ and a collection of simple closed curves $\alpha'$ such that $(\gamma',\mathfrak{a}')$ is slide equivalent to $(\gamma',\mathfrak{c})$, $\alpha'$ is slide equivalent to $\alpha$, and $\mathfrak{a}' \cap \alpha' = \emptyset$.
\item Take a collection of properly embedded arcs $\mathfrak{a}^*$ in $\Sigma$ such that $(\mathfrak{a}^*,\alpha)$ is slide equivalent to $(\mathfrak{a}',\alpha')$.
\end{enumerate}
The open book induced by $\calD$ has page $\Sigma_{\alpha}$, 
and the monodromy $\phi : \Sigma_{\alpha} \to \Sigma_{\alpha}$ is the unique (up to isotopy) diffeomorphism such that
\begin{equation*}
\phi(\mu_{\alpha}(\mathfrak{a})) = \mu_{\alpha}(\mathfrak{a}^*),
\end{equation*}
where the ordering of the arcs is preserved.
\end{theorem}

For a relative trisection diagram $(\Sigma;\alpha,\beta,\gamma)$, we call a collection $\mathfrak{a}$ of properly embedded arcs in $\Sigma$ an \textit{$\alpha$-cut system} if cutting $\Sigma_\alpha$ along $\mu_\alpha(\mathfrak{a})$ yields a disk.
Throughout this paper, the arcs in an $\alpha$-cut system are depicted in purple.
Theorem~\ref{thm:rtd-obd} gives an algorithm for describing the image of an $\alpha$-cut system under the induced monodromy.

\section{Splittable trisections}\label{sect:spl}

In this section, we introduce splittable trisections. We first define subdiagrams of trisection diagrams.

\begin{definition}
Let $\mathcal{D}=(\Sigma;\alpha,\beta,\gamma)$ be a (relative) trisection diagram.
Let $\Sigma'$ be a subsurface of $\Sigma$ such that
$\partial\Sigma'\cap(\alpha\cup\beta\cup\gamma)=\emptyset$.
We denote the diagram determined by the curve systems contained in $\Sigma'$ by $\mathcal{D}\cap\Sigma'$, namely,
\begin{equation*}
\mathcal{D}\cap\Sigma':=(\Sigma';\alpha\cap\Sigma',\beta\cap\Sigma',\gamma\cap\Sigma').
\end{equation*}
Here, $\alpha\cap\Sigma'$ denotes the collection of curves in $\alpha$ that are contained in $\Sigma'$, and similarly for $\beta$ and $\gamma$.
If $\mathcal{D}\cap\Sigma'$ is a relative trisection diagram, we call it a \emph{subdiagram} of $\mathcal{D}$.
%
\end{definition}

\begin{definition}\label{dfn:spl-rt}
A (relative) trisection diagram $\mathcal{D}=(\Sigma;\alpha,\beta,\gamma)$ is called \emph{splittable} if there exists a collection $\mathfrak{l}=\{\ell_1,\ldots,\ell_n\}$ of nontrivial simple closed curves on $\Sigma$ such that
\begin{itemize}
%
\item $\mathfrak{l}\cap(\alpha\cup\beta\cup\gamma)=\emptyset$.
\item Cutting $\Sigma$ along $\mathfrak{l}$ yields exactly two connected components, denoted by $\Sigma'$ and $\Sigma''$.
\item Both $\mathcal{D}\cap\Sigma'$ and $\mathcal{D}\cap\Sigma''$ are subdiagrams of $\calD$.
%
%
\end{itemize}
Any such collection $\mathfrak{l}$ is called a \emph{splitting system} for $\mathcal{D}$.
A relative trisection $\mathcal{T}$ is called \emph{splittable} if $\mathcal{T}$ admits a splittable trisection diagram.
\end{definition}

Although $\alpha$, $\beta$, $\gamma$, and $\mathfrak{l}$ are technically collections of curves, we treat them as subsets of $\Sigma$ in the definitions above. We will use this abuse of notation throughout the paper.

\begin{remark}
According to Definition~\ref{dfn:spl-rt}, every (relative) trisection diagram admits a trivial splitting, that is, it can be split into a small disk and its complement.
It may therefore be more appropriate to require that both $\mathcal{D}\cap\Sigma'$ and $\mathcal{D}\cap\Sigma''$ contain curve systems.
\end{remark}

\begin{figure}[!tbp]
\centering
\includegraphics[scale=0.23]{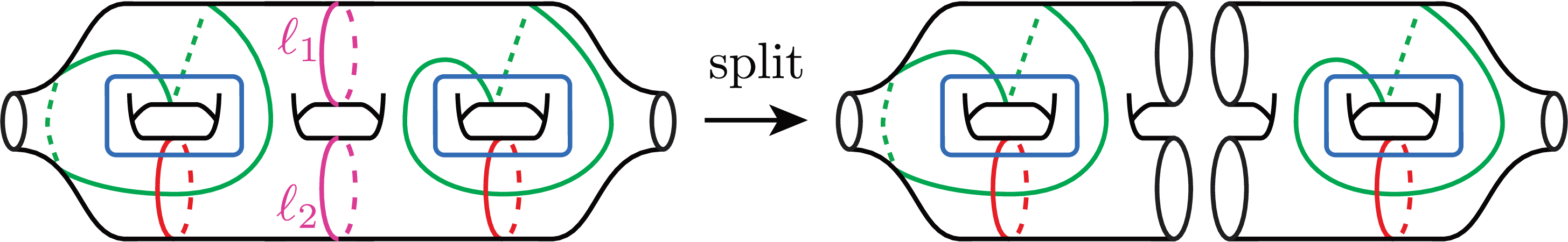}
\caption{A splittable relative trisection diagram.}
\label{fig:td-split}
\end{figure}

Figure~\ref{fig:td-split} shows a splittable relative trisection diagram with splitting system $\mathfrak{l}=\{\ell_1, \ell_2\}$.

We introduce a simple construction of a splittable trisection diagram.

\begin{lemma}\label{lem:glue-rtd}
Let $\mathcal{D}=(\Sigma;\alpha,\beta,\gamma)$ be a $(g,k;p,b)$-relative trisection diagram, and let $\mathfrak{l}=\{\ell_1,\ldots,\ell_n\}$ be an ordered collection of $n$ components of $\partial\Sigma$.
Similarly, let $\mathcal{D}'=(\Sigma';\alpha',\beta',\gamma')$ be a $(g',k';p',b')$-relative trisection diagram, and let $\mathfrak{l}'=\{\ell'_1,\ldots,\ell'_n\}$ be an ordered collection of $n$ components of $\partial\Sigma'$.
Let $\Sigma\cup_{\mathfrak{l}\sim\mathfrak{l}'}\Sigma'$ denote the surface obtained by identifying $\ell_i$ with $\ell'_i$ for each $i\in\{1,\ldots,n\}$, where the gluing is chosen to be compatible with the orientations of the two surfaces.
We define
\begin{equation*}
\mathcal{D}\cup_{\mathfrak{l}\sim\mathfrak{l}'}\mathcal{D}'
:=
\left(
\Sigma\cup_{\mathfrak{l}\sim\mathfrak{l}'}\Sigma';
\alpha\cup\alpha',
\beta\cup\beta',
\gamma\cup\gamma'
\right).
\end{equation*}
If the underlying surface $\Sigma\cup_{\mathfrak{l}\sim\mathfrak{l}'}\Sigma'$ has nonempty boundary, then $\mathcal{D}\cup_{\mathfrak{l}\sim\mathfrak{l}'}\mathcal{D}'$ is a relative trisection diagram of type
\begin{equation*}
(g+g'+n-1,k+k'-1;p+p'+n-1,b+b'-2n).
\end{equation*}
Moreover, $\mathcal{D}\cup_{\mathfrak{l}\sim\mathfrak{l}'}\mathcal{D}'$ is clearly splittable, with splitting system given by the image of $\mathfrak{l}$ and $\mathfrak{l}'$.
\end{lemma}

\begin{proof}
After identifying $\Sigma$ with $\Sigma_{g,b}$ and $\Sigma'$ with $\Sigma_{g',b'}$ by surface diffeomorphisms, we may assume that $\Sigma=\Sigma_{g,b}$ and $\Sigma'=\Sigma_{g',b'}$.
To prove that $\mathcal{D}\cup_{\mathfrak{l}\sim\mathfrak{l}'}\mathcal{D}'$ is a relative trisection diagram, it suffices to show that each of the diagrams
\begin{align*}
&\left(
\Sigma\cup_{\mathfrak{l}\sim\mathfrak{l}'}\Sigma';
\alpha\cup\alpha',\beta\cup\beta'
\right), \\
&\left(
\Sigma\cup_{\mathfrak{l}\sim\mathfrak{l}'}\Sigma';
\beta\cup\beta',\gamma\cup\gamma'
\right), \\
&\left(
\Sigma\cup_{\mathfrak{l}\sim\mathfrak{l}'}\Sigma';
\gamma\cup\gamma',\alpha\cup\alpha'
\right)
\end{align*}
is slide-diffeomorphic to the standard diagram shown in Figure~\ref{fig:td-gkpb}.
We verify this for the first diagram; the other two cases are analogous.
Since $\mathcal{D}$ and $\mathcal{D}'$ are relative trisection diagrams, there exist orientation-preserving diffeomorphisms $f:\Sigma\to\Sigma$ and $f':\Sigma'\to\Sigma'$ such that $(\Sigma;f(\alpha),f(\beta))$ and $(\Sigma';f'(\alpha'),f'(\beta'))$ are slide-equivalent to the standard diagrams.
We may assume that $f$ and $f'$ fix neighborhoods of the boundary pointwise.
Thus, for each $i\in\{1,\ldots,n\}$, the restrictions $f|_{\ell_i}$ and $f'|_{\ell'_i}$ agree under the identification $\ell_i=\ell'_i$.
Consequently, $f$ and $f'$ induce a well-defined diffeomorphism
\begin{equation*}
\widetilde{f}:
\Sigma\cup_{\mathfrak{l}\sim\mathfrak{l}'}\Sigma'
\longrightarrow
\Sigma\cup_{\mathfrak{l}\sim\mathfrak{l}'}\Sigma'.
\end{equation*}
satisfying $\widetilde{f}|_{\Sigma}=f$ and $\widetilde{f}|_{\Sigma'}=f'$.
It follows that
\begin{equation*}
\begin{aligned}
\widetilde{f}
\left(
\Sigma\cup_{\mathfrak{l}\sim\mathfrak{l}'}\Sigma';
\alpha\cup\alpha',\beta\cup\beta'
\right )
&=
\left(
\Sigma\cup_{\mathfrak{l}\sim\mathfrak{l}'}\Sigma';
\widetilde{f}(\alpha\cup\alpha'),
\widetilde{f}(\beta\cup\beta')
\right)\\
&=
\left(
\Sigma\cup_{\mathfrak{l}\sim\mathfrak{l}'}\Sigma';
\widetilde{f}(\alpha)\cup\widetilde{f}(\alpha'),
\widetilde{f}(\beta)\cup\widetilde{f}(\beta')
\right)\\
&=
\left(
\Sigma\cup_{\mathfrak{l}\sim\mathfrak{l}'}\Sigma';
f(\alpha)\cup f'(\alpha'),
f(\beta)\cup f'(\beta')
\right).
\end{aligned}
\end{equation*}
Therefore, the resulting diagram is obtained by gluing $(\Sigma;f(\alpha),f(\beta))$ and $(\Sigma';f'(\alpha'),f'(\beta'))$ via the identifications $\ell_i=\ell'_i$ for each $i\in\{1,\ldots,n\}$.
Each of the two pieces separated by the gluing circles can be put into standard form by a sequence of slides. See Figure~\ref{fig:td-prf-glue}.
\begin{figure}[!tbp]
\centering
\includegraphics[scale=0.23]{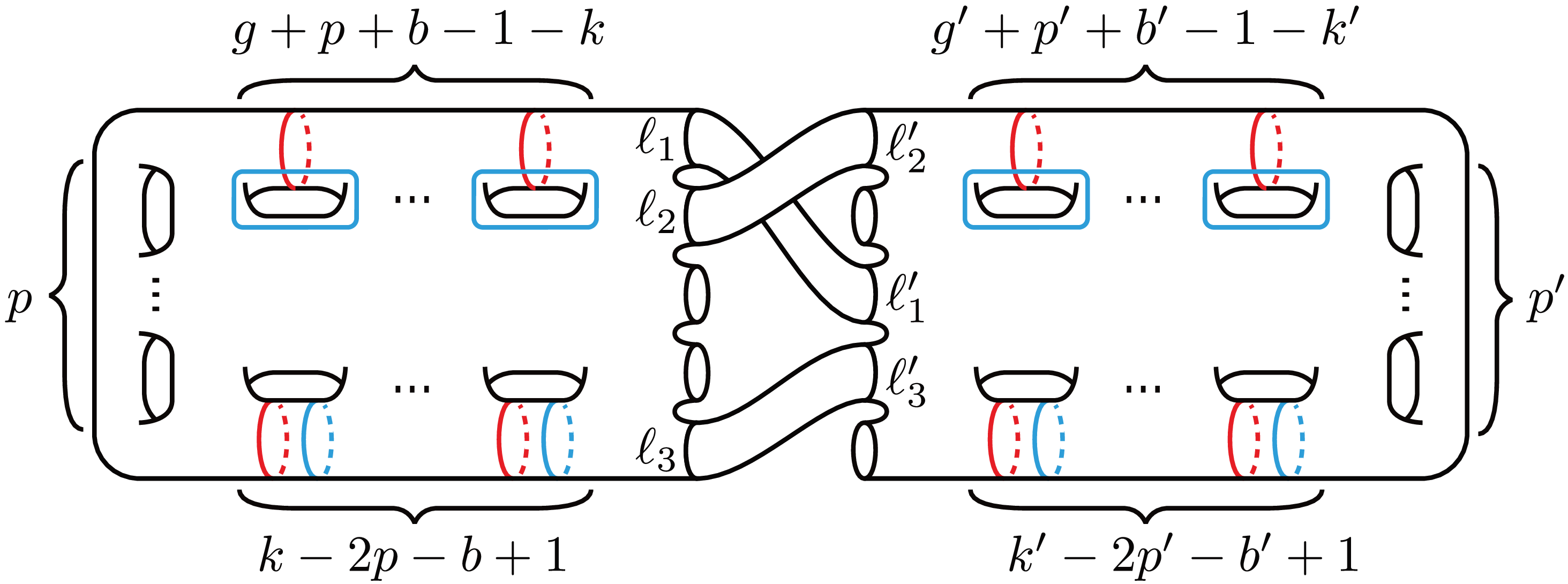}
\caption{The diagram $\left( \Sigma\cup_{\mathfrak{l}\sim\mathfrak{l}'}\Sigma'; f(\alpha)\cup f'(\alpha'), f(\beta)\cup f'(\beta') \right)$ after putting each piece into standard form by slides. The case $b=b'=5$ and $n=3$.}
\label{fig:td-prf-glue}
\end{figure}
We see that the resulting diagram is diffeomorphic to the standard diagram of type
\begin{equation*}
(g+g'+n-1,k+k'-1;p+p'+n-1,b+b'-2n).
\end{equation*}
Indeed, the underlying surface $\Sigma\cup_{\mathfrak{l}\sim\mathfrak{l}'}\Sigma'$ has genus $g+g'+n-1$ and $b+b'-2n$ boundary components.
The numbers of components in each curve system, dual pairs, and parallel pairs are given, respectively, by
\begin{align*}
& (g+g'+n-1)-(p+p'+n-1), \\
& (g+g'+n-1)+(p+p'+n-1)+(b+b'-2n)-1-(k+k'-1), \\
& (k+k'-1)-2(p+p'+n-1)-(b+b'-2n)+1.
\end{align*}
\end{proof}

\begin{figure}[!tbp]
\centering
\includegraphics[scale=0.23]{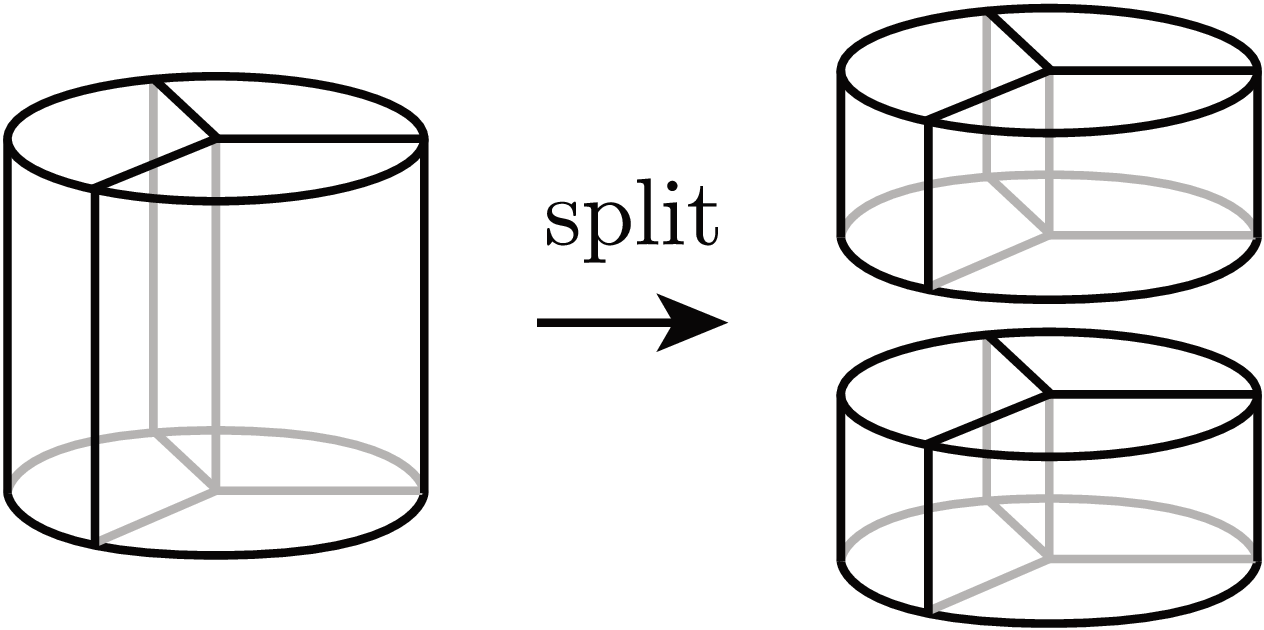}
\caption{A schematic picture of a splittable relative trisection.}
\label{fig:split}
\end{figure}

Let $X=X_1\cup X_2\cup X_3$ and $X'=X'_1\cup X'_2\cup X'_3$ be the relative trisections represented by $\mathcal{D}$ and $\mathcal{D}'$, respectively.
We construct a relative trisection represented by $\mathcal{D}\cup_{\mathfrak{l}\sim\mathfrak{l}'}\mathcal{D}'$ as follows.
By the definition of a relative trisection, each $X_i$ is diffeomorphic to $Z_k = (\Sigma_{p,b}\times\Delta)\natural V_n$.
Similarly, each $X'_i$ is diffeomorphic to $Z_k' = (\Sigma_{p',b'}\times\Delta)\natural V_{n'}$.
Let $X''_i$ be the $4$-manifold obtained by gluing $X_i$ and $X'_i$ via the identification of $\mathfrak{l}\times\Delta$ with $\mathfrak{l}'\times\Delta$.
Then
\begin{equation*}
X''_i
\cong
\left(
(\Sigma_{p,b}\cup_{\mathfrak{l}\sim\mathfrak{l}'}\Sigma_{p',b'})
\times\Delta
\right)
\natural (V_n\natural V_{n'})
\cong Z_{k+k'-1}
\end{equation*}
It follows that $X''_1\cup X''_2\cup X''_3$ is a relative trisection, which is precisely the gluing of $X=X_1\cup X_2\cup X_3$ and $X'=X'_1\cup X'_2\cup X'_3$. See Figure~\ref{fig:split}.
Moreover, by construction, $\mathcal{D}\cup_{\mathfrak{l}\sim\mathfrak{l}'}\mathcal{D}'$ is a relative trisection diagram for $X''_1\cup X''_2\cup X''_3$.

In the construction of a splittable relative trisection in Lemma~\ref{lem:glue-rtd}, we show that the induced open book is also splittable, particularly when gluing along a single boundary component.

\begin{lemma}\label{lem:srtd-obd}
Let $\mathcal{D}=(\Sigma;\alpha,\beta,\gamma)$ be a $(g,k;p,b)$-relative trisection diagram with $b\geq2$, and let $\ell$ be a component of $\partial\Sigma$.
Similarly, let $\mathcal{D}'=(\Sigma';\alpha',\beta',\gamma')$ be a $(g',k';p',b')$-relative trisection diagram with $b'\geq2$, and let $\ell'$ be a component of $\partial\Sigma'$.
Let $\mathcal{D}\cup_{{\ell}\sim\ell'}\mathcal{D}'$ be the $(g+g', k+k'-1; p+p', b+b'-2)$ relative trisection diagram defined in Lemma~\ref{lem:glue-rtd}.
Let $\phi:\Sigma_\alpha\to\Sigma_\alpha$ and $\phi':\Sigma'_{\alpha'}\to\Sigma'_{\alpha'}$ be the monodromies induced by $\mathcal{D}$ and $\mathcal{D}'$, respectively.
Let
\begin{equation*}
\Phi:
\Sigma_\alpha\cup_{\ell\sim\ell'}\Sigma'_{\alpha'}
\longrightarrow
\Sigma_\alpha\cup_{\ell\sim\ell'}\Sigma'_{\alpha'}
\end{equation*}
be the map defined by $\Phi(x)=\phi(x)$ if $x\in\Sigma_\alpha$ and $\Phi(x)=\phi'(x)$ if $x\in\Sigma'_{\alpha'}$.
Then $\Phi$ is the monodromy induced by $\mathcal{D}\cup_{\ell\sim\ell'}\mathcal{D}'$.
\end{lemma}

\begin{proof}
Let $\frak{a}=\{a_1,\ldots,a_{2p+b-1}\}$ be a cut system for $\mathcal{D}$.
Namely, cutting the surgered surface $\Sigma_{\alpha}$ along the arcs in $\frak{a}$ yields a disk.
We may assume that $\ell$ intersects $\frak{a}$ only in $a_1$, i.e., $\ell\cap\frak{a}=\ell\cap a_1=\{ \mathrm{pt} \}$.
Applying the algorithm of Theorem~\ref{thm:rtd-obd} to $\mathcal{D}$, we obtain can draw the image
\begin{equation*}
\phi(\frak{a})=\{\phi(a_1),\ldots,\phi(a_{2p+b-1})\}
\end{equation*}
on $\Sigma_{\alpha}$ up to isotopy.
We proceed similarly for $\mathcal{D}'$, assuming that the diagrams are glued along $\ell=\ell'$ so that $\ell\cap a_1=\ell'\cap a'_1$.
Note that the monodromies fix neighborhoods of the boundaries.

Define a collection $\frak{a}\cup_{\ell\sim\ell'}\frak{a}'$ of $2(p+p')+(b+b'-2)-1$ arcs on $\Sigma\cup_{\ell\sim\ell'}\Sigma'$ by
\begin{equation*}
\frak{a}\cup_{\ell\sim\ell'}\frak{a}' :=
\{ a_1\cup_{\ell\sim\ell'}a'_1,
a_2,\ldots,a_{2p+b-1}, a'_2,\ldots,a'_{2p'+b'-1}\}.
\end{equation*}
Here, $a_1\cup_{\ell\sim\ell'}a'_1$ is the arc obtained by gluing $a_1$ and $a'_1$ along the identification of the points $\ell\cap a_1=\ell'\cap a'_1$.
Then $\frak{a}\cup_{\ell\sim\ell'}\frak{a}'$ is a cut system for $\mathcal{D}\cup_{\ell\sim\ell'}\mathcal{D}'$,
that is, cutting $\Sigma_\alpha\cup_{\ell\sim\ell'}\Sigma'_{\alpha'}$ along these arcs yields a disk.
By applying the algorithm of Theorem~\ref{thm:rtd-obd} to $\mathcal{D}\cup_{\ell\sim\ell'}\mathcal{D}'$, we can draw the image of $\frak{a}\cup_{\ell\sim\ell'}\frak{a}'$ under the induced monodromy on $\Sigma_\alpha\cup_{\ell\sim\ell'}\Sigma'_{\alpha'}$ as
\begin{equation*}
\begin{aligned}
\{ \phi(a_1)\cup_{\ell\sim\ell'}\phi'(a'_1),
\phi(a_2),\ldots,\phi(a_{2p+b-1}), \phi'(a'_2),\ldots, \phi'(a'_{2p+b-1}) \}
\end{aligned}
\end{equation*}
up to isotopy.
Indeed, since the arcs $a_2,\ldots,a_{2p+b-1}$ lie entirely in $\mathcal{D}$ and $a'_2,\ldots,a'_{2p'+b'-1}$ lie entirely in $\mathcal{D}'$, the algorithm can be applied to these arcs independently, exactly as in the original diagrams.
For the remaining arc $a_1\cup_{\ell\sim\ell'}a'_1$, we first apply the algorithm to the portion lying in $\mathcal{D}$, and then apply it to the portion lying in $\mathcal{D}'$.
This yields the arc $\phi(a_1)\cup_{\ell\sim\ell'}\phi'(a'_1)$.
The resulting image agrees with the image of $\frak{a}\cup \frak{a}'$ under $\Phi$, as desired.
\end{proof}

\section{Trisecting contractible 4-manifolds}\label{sect:tris-cork}

We give relative trisections of the contractible $4$-manifolds $C_n$ introduced in Subsection~\ref{subsec:corks}.

\begin{figure}[!tbp]
\centering
\includegraphics[scale=0.23]{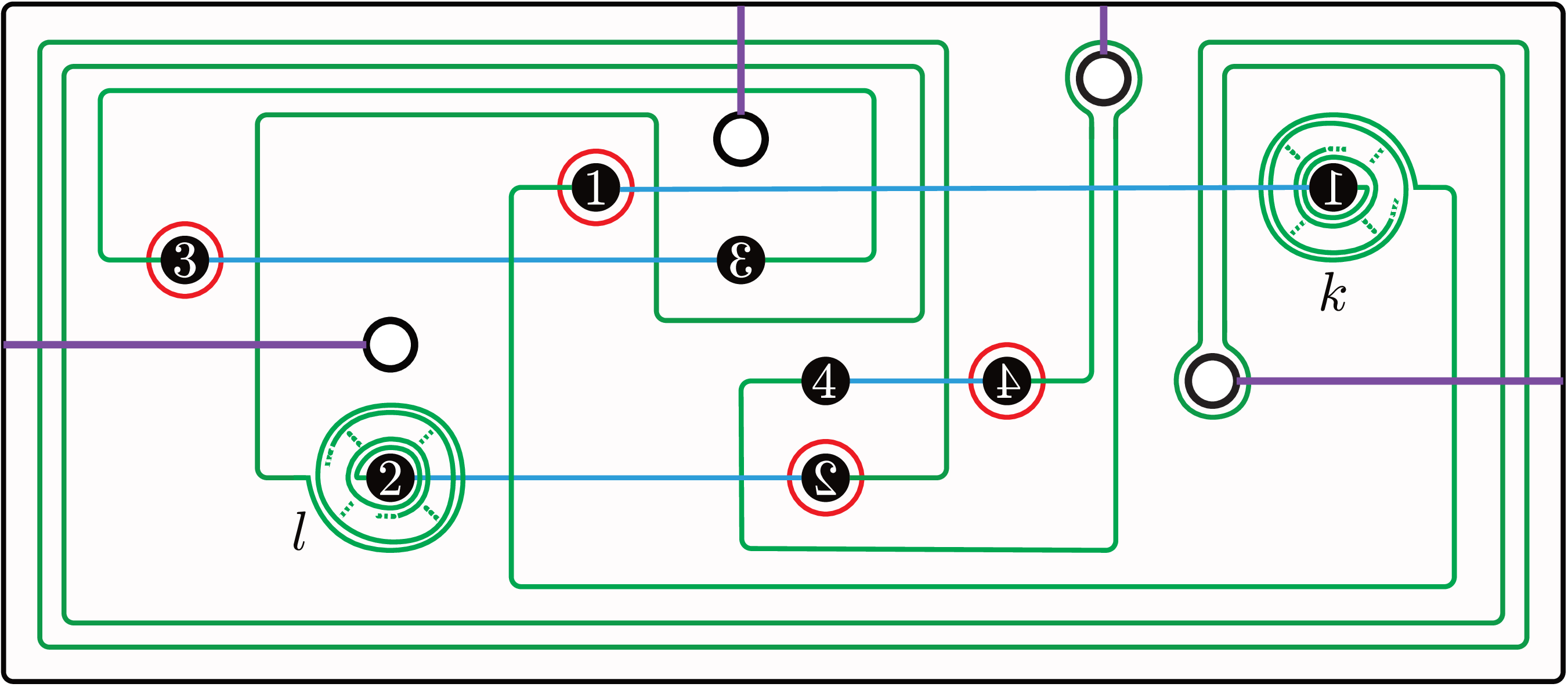}
\caption{The diagram $\mathcal{D}_{k,l}$. The purple arcs form a cut system.}
\label{fig:td-Cn}
\end{figure}

\begin{definition}
For integers $k$ and $l$, let $\mathcal{D}_{k,l}:=(\Sigma;\alpha,\beta,\gamma)$ be the diagram shown in Figure~\ref{fig:td-Cn}, where the $\gamma$-curves are obtained by applying $k$ positive Dehn twists along the boundary of the black disk labeled ``$1$'' and $l$ positive Dehn twists along the boundary of the black disk labeled ``$2$''.
If either $k$ or $l$ is negative, the corresponding Dehn twists are negative.
\end{definition}

\begin{proposition}\label{prop:rtd-Cn}
For any integers $k$ and $l$, $\mathcal{D}_{k,l}$ is a $(4,4;0,5)$-relative trisection diagram of $C_{k+l+2}$.
\end{proposition}

\begin{figure}[!tbp]
\centering
\includegraphics[scale=0.23]{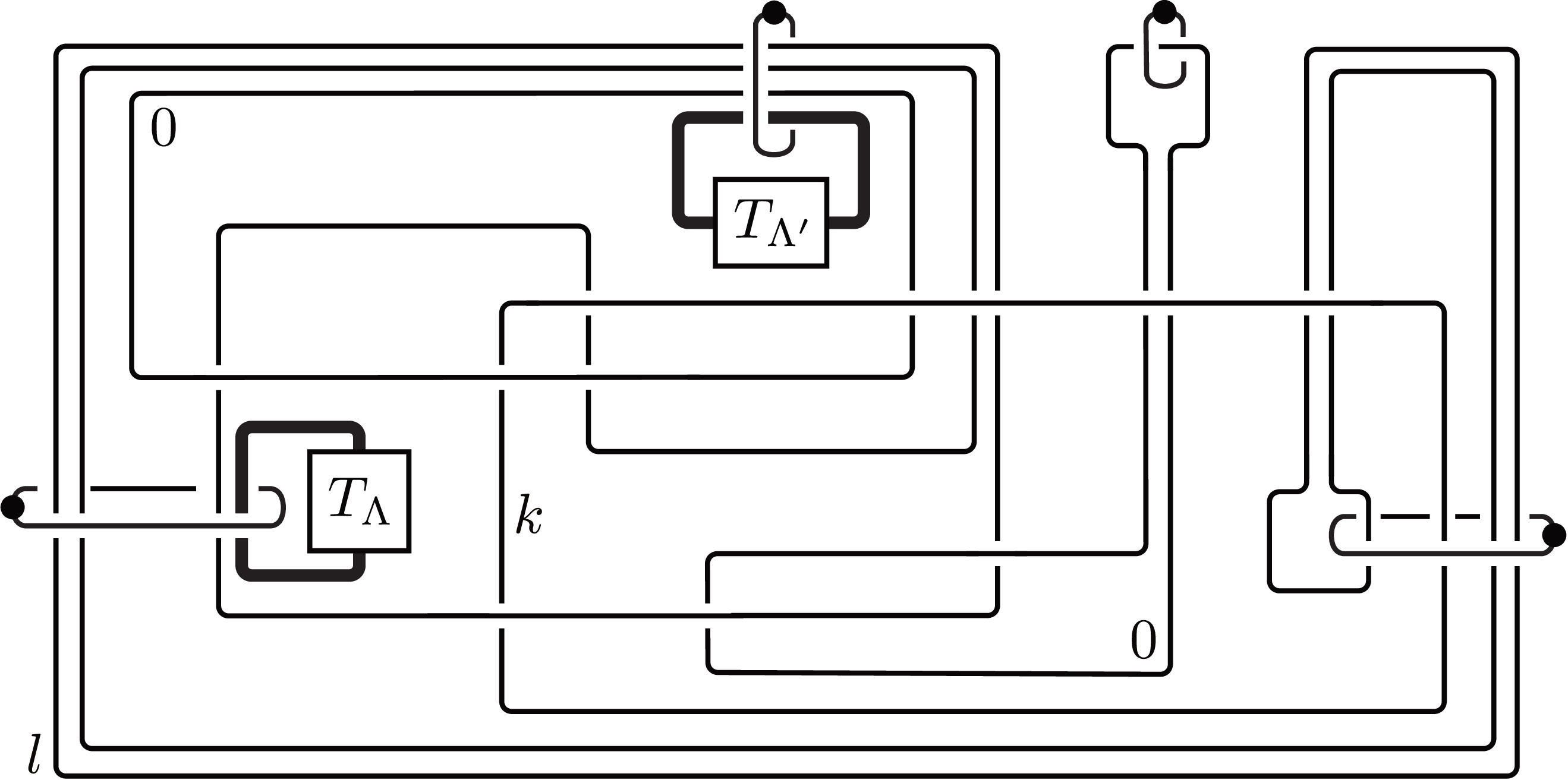}
\caption{The handlebody diagram induced by $\calD_{k,l}$ with tangles $T_{\Lambda}$ and $T_{\Lambda'}$.}
\label{fig:kd-induced}
\end{figure}

\begin{proof}
First, we prove that $\mathcal{D}_{k,l}$ is a $(4,4;0,5)$-relative trisection diagram.
It suffices to check that each of the triples $(\Sigma;\alpha,\beta)$, $(\Sigma;\beta,\gamma)$, and $(\Sigma;\gamma,\alpha)$ is slide-diffeomorphic to the standard diagram shown in Figure~\ref{fig:td-gkpb}.
Both $(\Sigma;\alpha,\beta)$ and $(\Sigma;\gamma,\alpha)$ can be put into standard form by diffeomorphisms of $\Sigma$, since every curve in each diagram has a unique dual curve.

For the remaining triple $(\Sigma;\beta,\gamma)$, we prove that it is slide-diffeomorphic to the standard diagram using the three operations shown in Figure~\ref{fig:5opes}.
\begin{figure}[!tbp]
\centering
\includegraphics[scale=0.23]{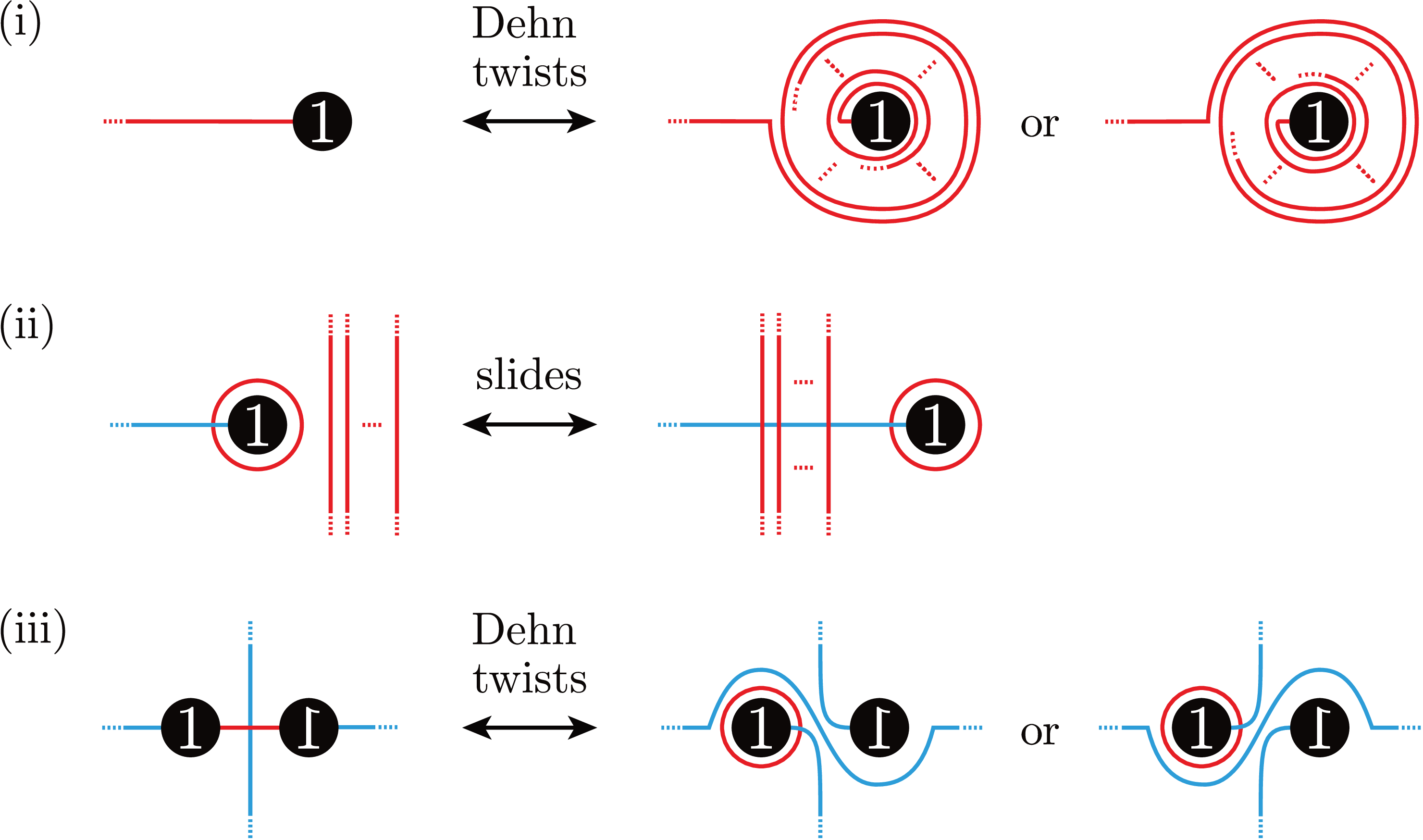}
\caption{Operations obtained by surface diffeomorphisms and slides of curves.}
\label{fig:5opes}
\end{figure}
See Section~3 of \cite{Tak25b} for details of there operations.
Figures~\ref{fig:td-Sbc1} and~\ref{fig:td-Sbc2} show a sequence of moves that modifies $(\Sigma;\beta,\gamma)$ into the standard diagram.
%
%
In the first diagram of Figure~\ref{fig:td-Sbc1}, we apply the operation~(iii) to the highlighted region.
In the second diagram, we drag the black disk labeled ``4'' and enclosed by a $\beta$-circle along the $\gamma$-curve as indicated.
The black disk can pass through the horizontal $\beta$-curves by applying the operation~(ii).
Note that the operation~(i) allows us to disregard the number of rotations of the $\gamma$-curve around the black disk.
In the third diagram, we apply the operation~(iii) to the highlighted region.
In the fourth diagram, we drag the black disk labeled ``3'' and enclosed by a $\beta$-circle along the $\gamma$-curve as indicated.
We then move to Figure~\ref{fig:td-Sbc2}. In the first diagram, we drag the hole to the indicated position.
In the second diagram, we drag the black disk labeled ``1'' and enclosed by a $\gamma$-circle along the $\beta$-curve as indicated.
In the third diagram, we drag the hole.
In the final diagram of Figure~\ref{fig:td-Sbc2}, every curve has a unique dual curve. Hence, this diagram is diffeomorphic to the standard diagram.
\begin{figure}[!tbp]
\centering
\includegraphics[scale=0.23]{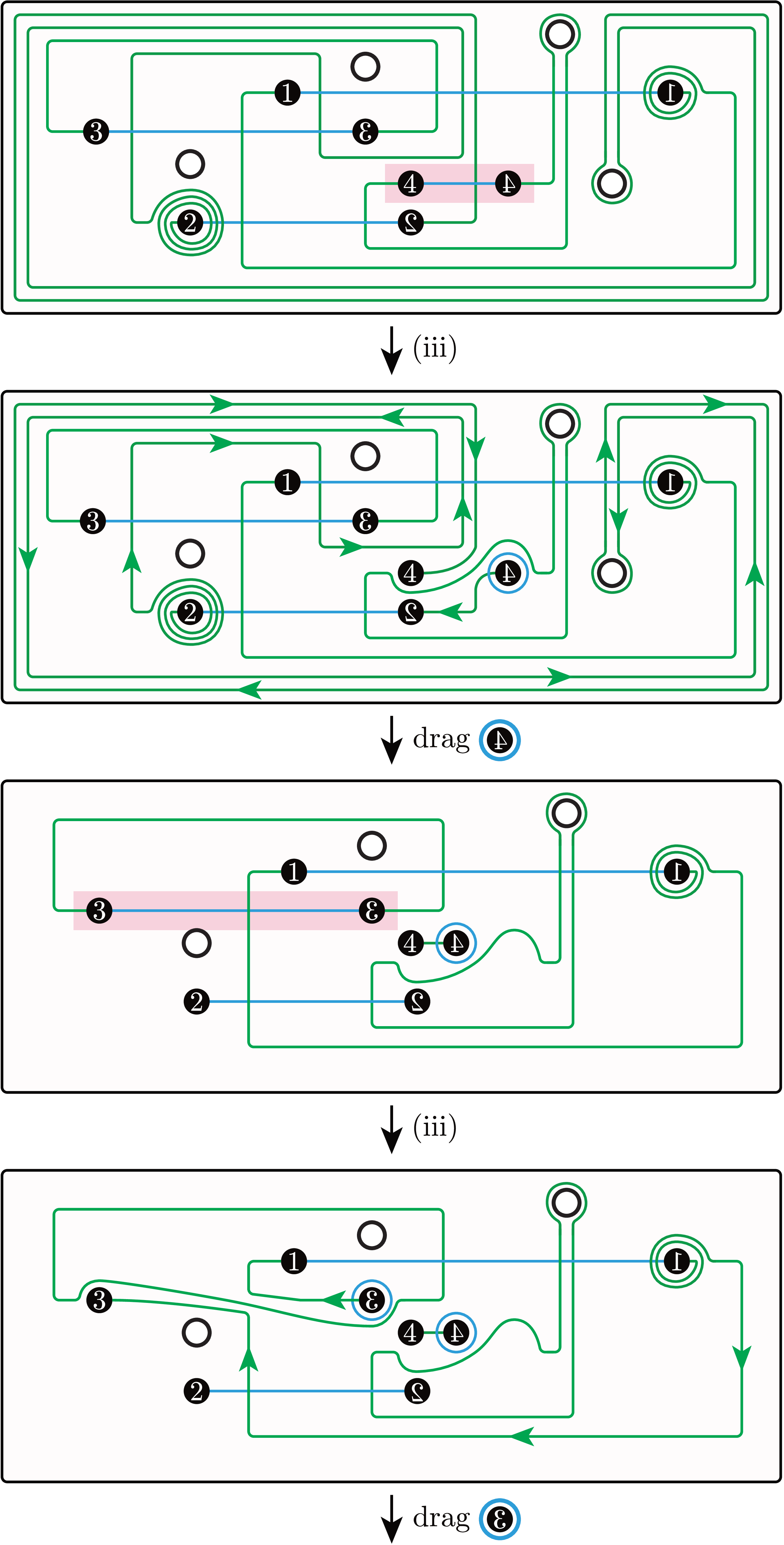}
\caption{Surface diffeomorphisms and slides of curves modifying $(\Sigma;\beta, \gamma)$ into the standard diagram. ($k=2$ and $l=3$)}
\label{fig:td-Sbc1}
\end{figure}
\begin{figure}[!tbp]
\centering
\includegraphics[scale=0.23]{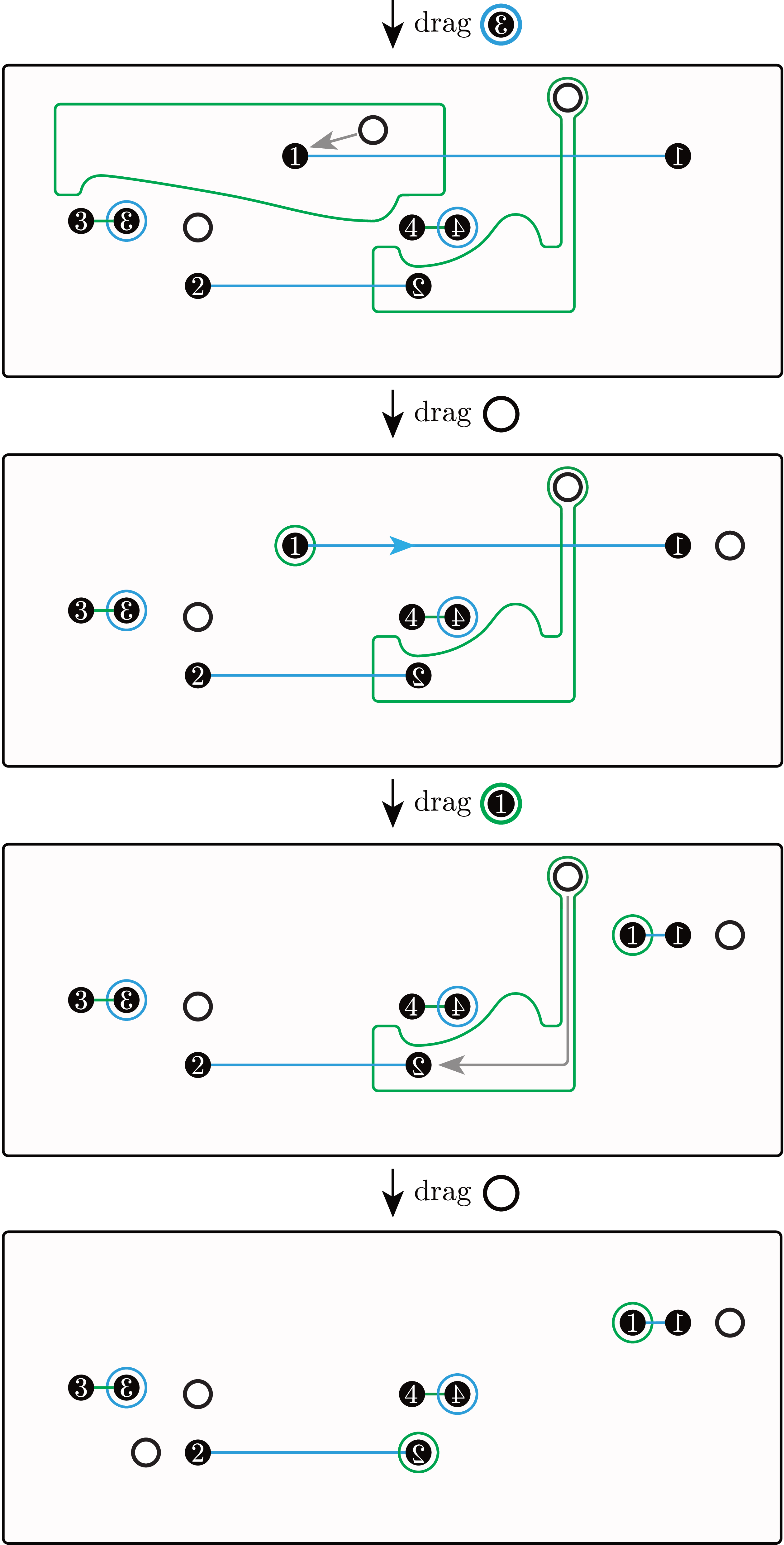}
\caption{Continuation of Figure~\ref{fig:td-Sbc1}.}
\label{fig:td-Sbc2}
\end{figure}

Next, we prove that the $4$-manifold represented by $\mathcal{D}_{k,l}$ is diffeomorphic to $C_{k+l+2}$.
To this end, we apply the algorithm of Kim--Miller~\cite[Section~2.2]{KimMil20}, which converts a relative trisection diagram into a handlebody diagram.
Applying this algorithm to $\mathcal{D}_{k,l}$ yields the handlebody diagram shown in Figure~\ref{fig:kd-induced} with the tangles $T_\Lambda$ and $T_\Lambda'$ removed.
Note that each purple cut arc is converted into a dotted circle, and each $\gamma$-curve is converted into the attaching circle of a $2$-handle with the surface framing.
By Lemma~\ref{lem:kcalc} below, the resulting diagram represents $C_{k+l+2}$.
\end{proof}

\begin{lemma}\label{lem:kcalc}
The $4$-manifold obtained by removing both tangles $T_\Lambda$ and $T_\Lambda'$ from the handlebody diagram in Figure~\ref{fig:kd-induced} is diffeomorphic to $C_{k+l+2}$.
The $4$-manifolds obtained by removing either $T_\Lambda$ or $T_\Lambda'$ from the handlebody diagram in Figure~\ref{fig:kd-induced} are diffeomorphic to $Y_1(\Lambda)$ and $Y_2(\Lambda')$, respectively.
\end{lemma}

\begin{figure}[!tbp]
\centering
\includegraphics[scale=0.23]{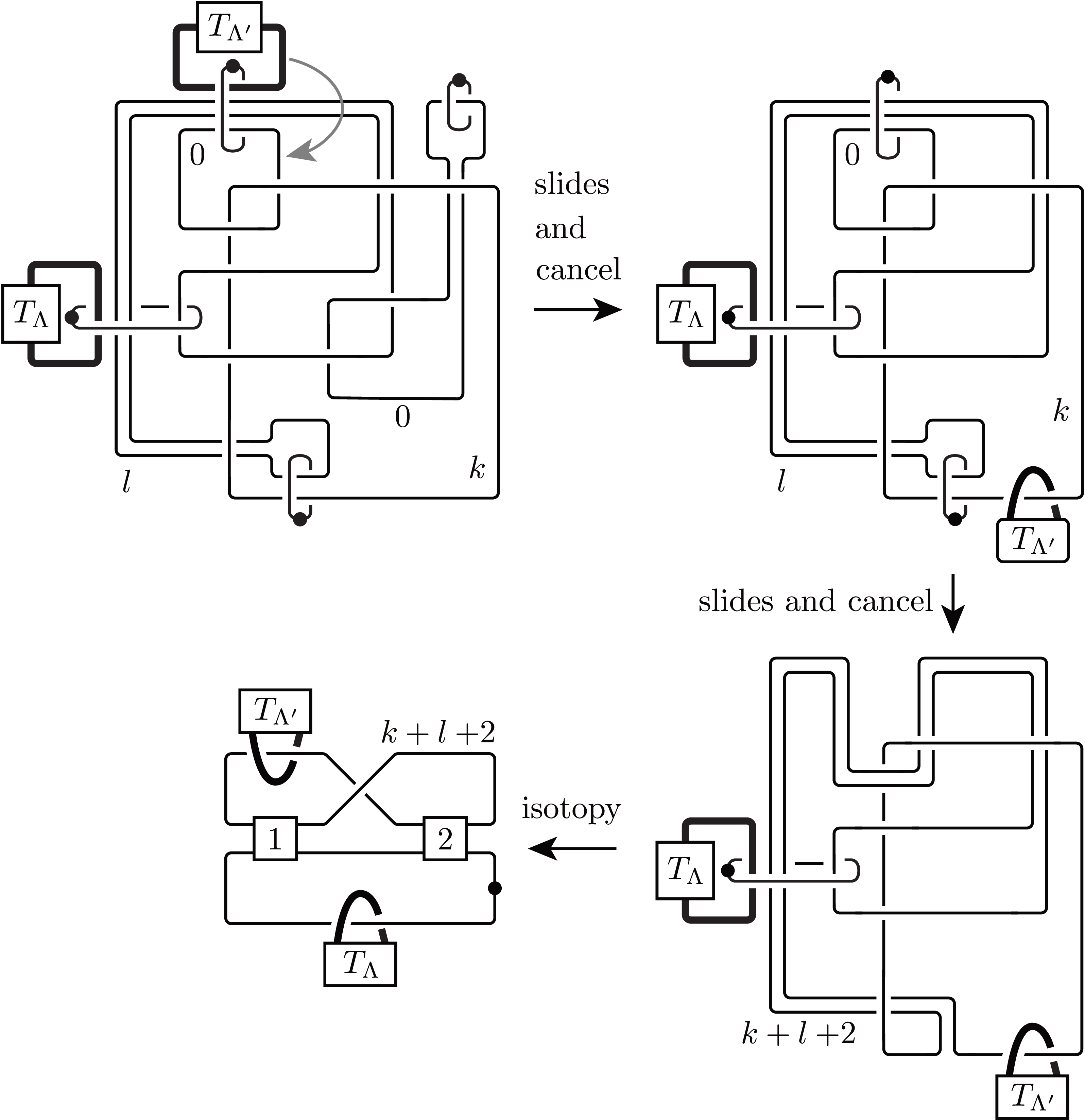}
\caption{Handlemoves of $C_n$ with tangles.}
\label{fig:kd-calc}
\end{figure}

\begin{proof}
After an isotopy of the diagram in Figure~\ref{fig:kd-induced}, we obtain the first diagram in Figure~\ref{fig:kd-calc}.
The proof is completed by the subsequent handle moves shown in Figure~\ref{fig:kd-calc}.
\end{proof}

\begin{remark}
By Proposition~\ref{prop:rtd-Cn}, the contractible 4-manifold $C_n$ admits a relative trisection of genus 4, but this is not the minimal genus.
In \cite{Tak25a}, a genus-$3$ relative trisection of $C_n$ is constructed, and this has minimal genus at least when $n \notin \{2,3,4,5\}$.
When $n \in \{ 2,3,4,5 \}$, $C_n$ could potentially admit a genus-2 relative trisection.
\end{remark}

\begin{remark}
The $(4,4;0,5)$-relative trisection diagram $\mathcal{D}_{k,l}$ induces an open book decomposition of the boundary homology $3$-sphere with page $\Sigma_{0,5}$.
Although the monodromy induced by $\mathcal{D}_{k,l}$ can be determined by following the procedure described in Theorem~\ref{thm:rtd-obd}, we do not specify it here, as carrying out this procedure is rather involved.
\end{remark}

\section{Constructions of splittable trisections for exotic pairs}\label{sect:tris-exotic}

In this section, we construct splittable relative trisection diagrams for exotic 4-manifolds related by twisting along the Akbulut cork.

\subsection{Producing homeomorphic pairs}\label{sect:homeo}

\begin{figure}[!tbp]
\centering
\includegraphics[scale=0.23]{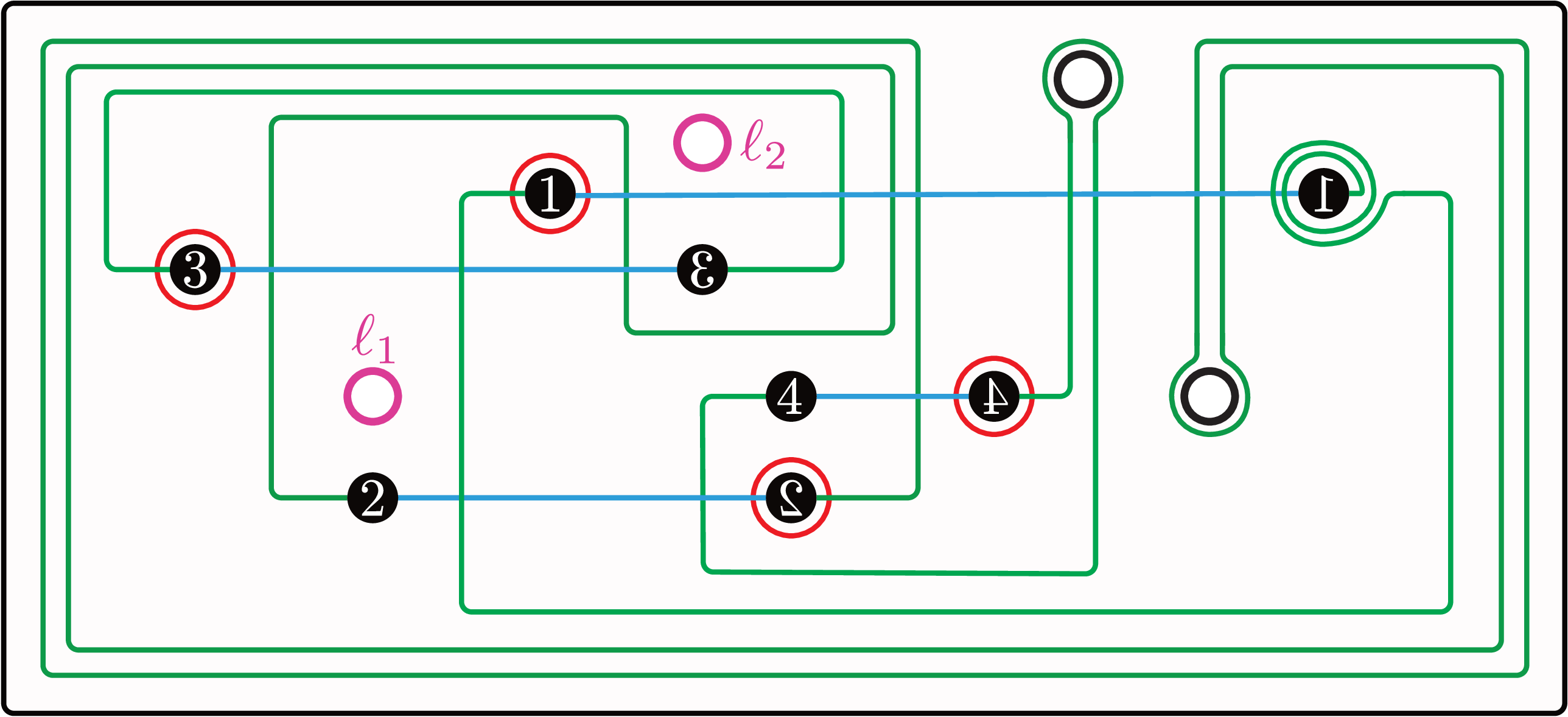}
\caption{$(4,4;0,5)$-relative trisection diagram $\mathcal{D}_A$ of the Akbulut cork $C_0$ with distinguished boundary components $\ell_1$ and $\ell_2$.}
\label{fig:td-Akb}
\end{figure}

In this subsection, we prove Theorem~\ref{main:td-homeo}.
Let $\mathcal{D}_A=\mathcal{D}_{-2,0}$, that is, let $\mathcal{D}_A$ be the $(4,4;0,5)$-relative trisection diagram of the Akbulut cork $C_0$ shown in Figure~\ref{fig:td-Akb}.
By Proposition~\ref{prop:rtd-Cn}, the diagram $\mathcal{D}_{k,l}$ represents $C_{k+l+2}$. Thus, we may use any other diagram $\mathcal{D}_{k,l}$ satisfying $k+l+2=0$ in place of $\mathcal{D}_{-2,0}$.
Let $\ell_1$ and $\ell_2$ be the distinguished boundary components shown in Figure~\ref{fig:td-Akb}.
The genus-$4$ pieces in Figures~\ref{fig:td-exoticY1} and~\ref{fig:td-exoticY2} are obtained from the diagram in Figure~\ref{fig:td-Akb} by applying a diffeomorphism of the surface.

\begin{lemma}\label{lem:glue-Akb-Lambda}
Let $\mathcal{D}$ be a relative trisection diagram with a distinguished boundary component $\ell$. Suppose that capping off $\ell$ with a disk yields a relative trisection diagram inducing a $2$-handlebody $\Lambda$.
For each $i\in\{1,2\}$, the diagram $\mathcal{D}_A\cup_{\ell_i\sim\ell}\mathcal{D}$, obtained by identifying $\ell_i$ with $\ell$, represents the $4$-manifold $Y_i(\Lambda)$.
\end{lemma}

\begin{figure}[!tbp]
\centering
\includegraphics[scale=0.23]{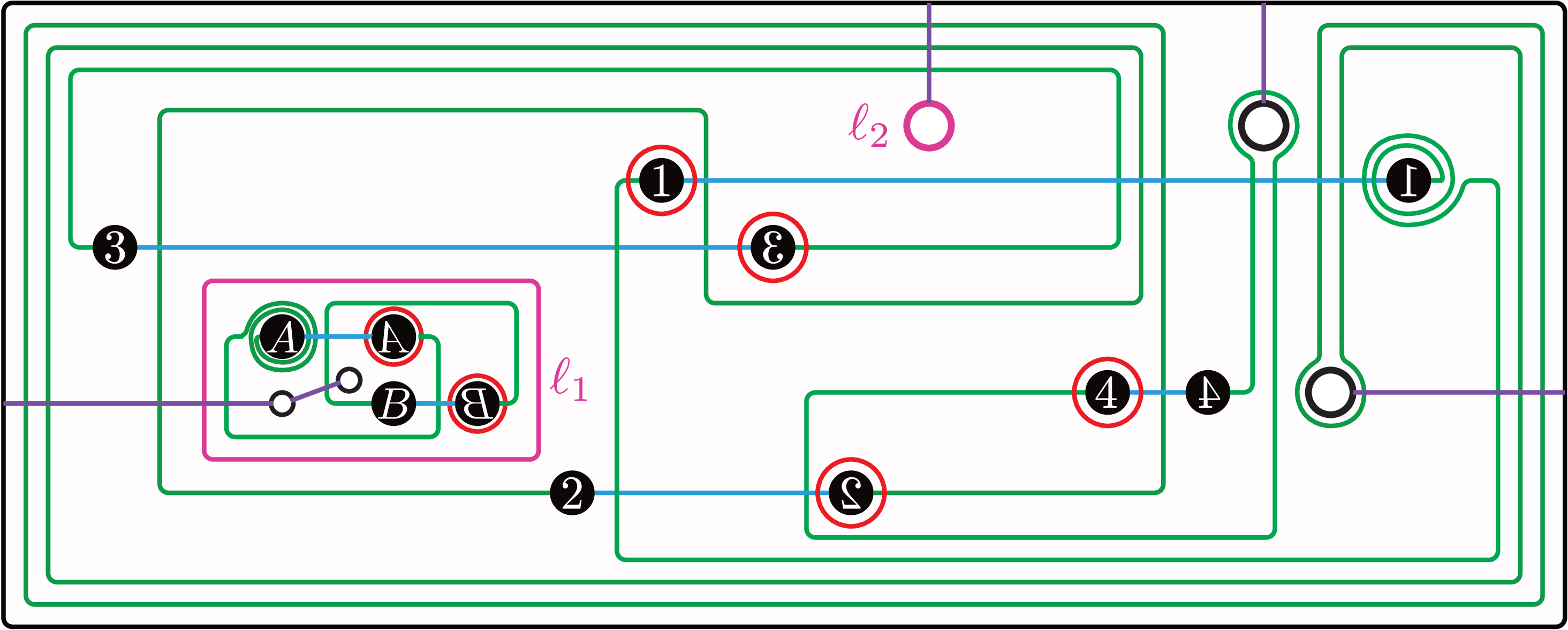}
\caption{$\mathcal{D}_A\cup_{\ell_1\sim\ell}\mathcal{D}$}
\label{fig:td-Ac+U21}
\end{figure}

\begin{figure}[!tbp]
\centering
\includegraphics[scale=0.23]{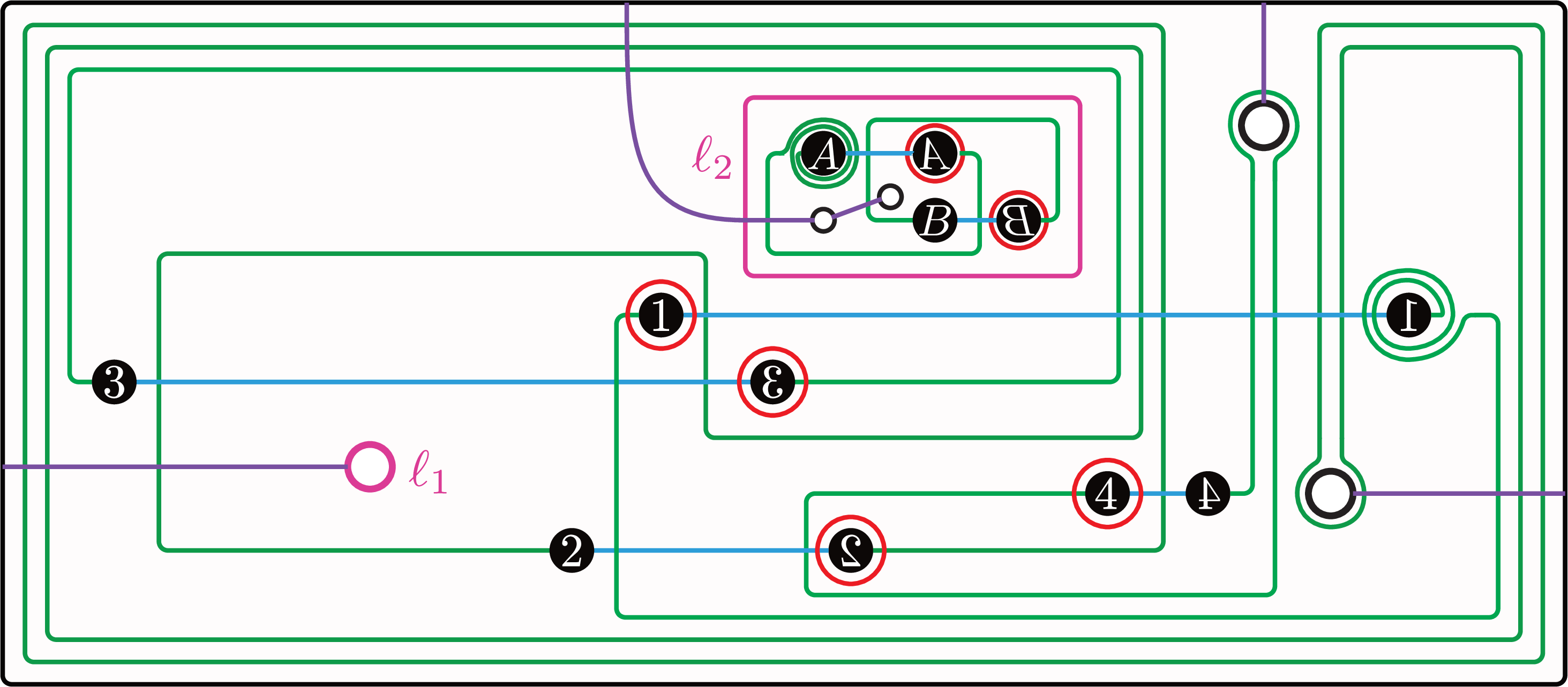}
\caption{$\mathcal{D}_A\cup_{\ell_2\sim\ell}\mathcal{D}$}
\label{fig:td-Ac+U22}
\end{figure}

Figures~\ref{fig:td-Ac+U21} and~\ref{fig:td-Ac+U22} show examples of $\mathcal{D}_A\cup_{\ell_i\sim\ell}\mathcal{D}$, where $\mathcal{D}$ in each figure is the diagram shown in Figure~\ref{fig:td-exm}.

\begin{figure}[!tbp]
\centering
\includegraphics[scale=0.23]{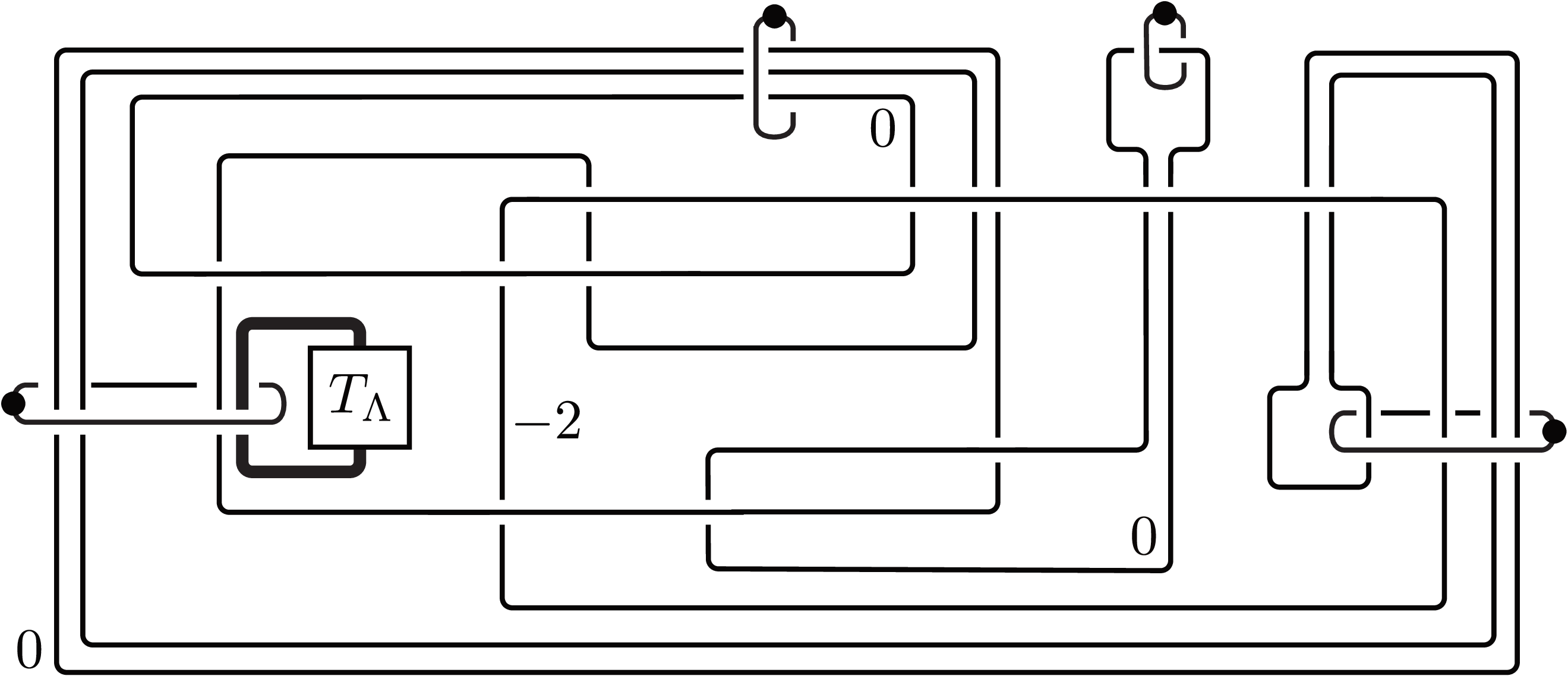}
\caption{The handlebody diagram induced by $\mathcal{D}_A\cup_{\ell_1\sim\ell}\mathcal{D}$.}
\label{fig:kd-Ac+L1}
\end{figure}

\begin{figure}[!tbp]
\centering
\includegraphics[scale=0.23]{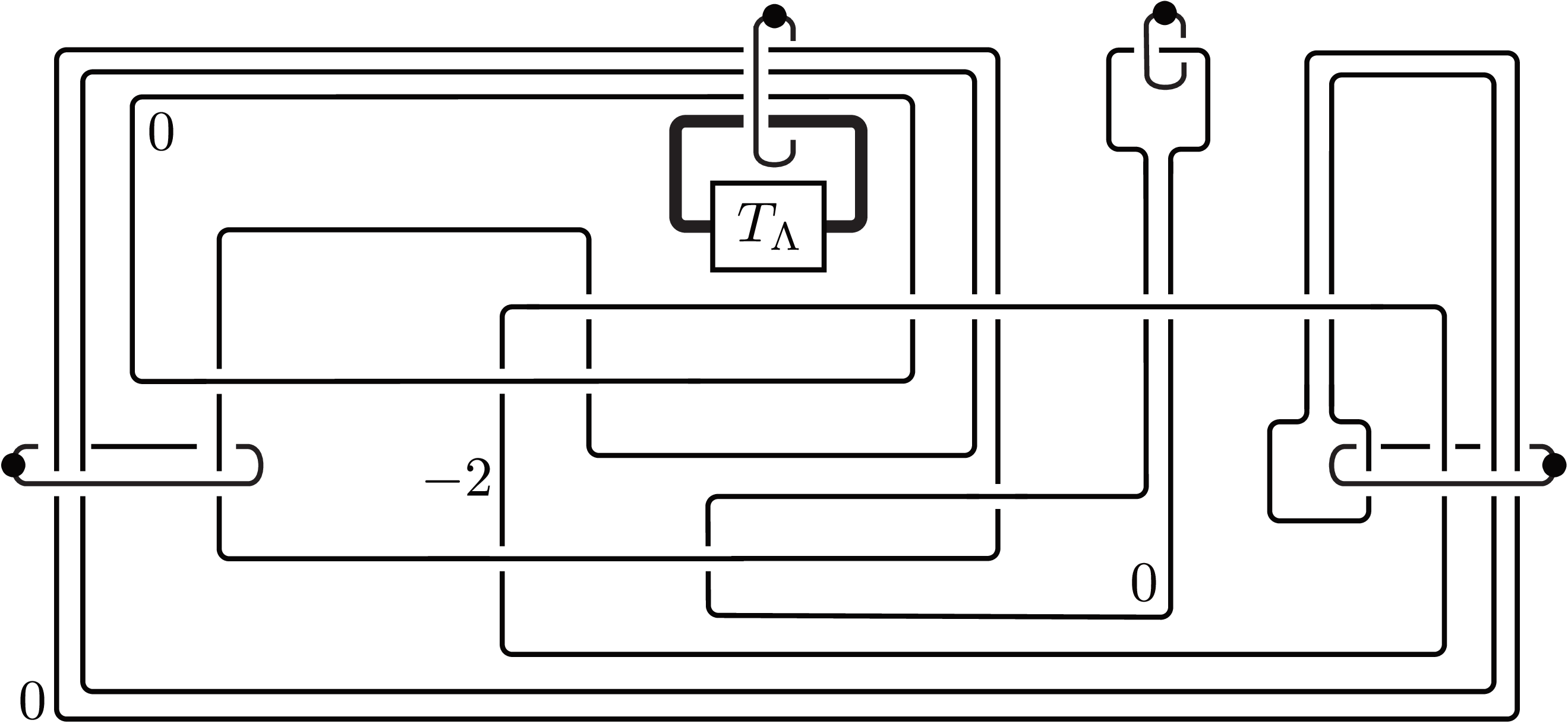}
\caption{The handlebody diagram induced by $\mathcal{D}_A\cup_{\ell_2\sim\ell}\mathcal{D}$.}
\label{fig:kd-Ac+L2}
\end{figure}

\begin{proof}[Proof of Lemma~\ref{lem:glue-Akb-Lambda}]
By Lemma~\ref{lem:glue-rtd}, any diagram obtained by such a gluing is a relative trisection diagram.
By applying the Kim--Miller's algorithm, we see that the resulting relative trisection diagrams induce the handlebody diagrams shown in Figures~\ref{fig:kd-Ac+L1} and~\ref{fig:kd-Ac+L2}, respectively.
The result now follows from Lemma~\ref{lem:kcalc}.
\end{proof}

\begin{proof}[Proof of Theorem~\ref{main:td-homeo}]
By Lemma~\ref{lem:glue-Akb-Lambda}, the 4-manifolds obtained by these gluings are of the form $Y_1(\Lambda)$ and $Y_2(\Lambda)$.
The 4-manifolds $Y_1(\Lambda)$ and $Y_2(\Lambda)$ are homeomorphic, since they are related by a cork twist.
\end{proof}

\subsection{Trisecting 2-handles}

The exotic $4$-manifolds studied in this paper are obtained by attaching a single $2$-handle to the Akbulut cork $C_0$.
In this subsection, we give relative trisections of $D^4$ that describe the $2$-handles to be attached to $C_0$.
In the corresponding handlebody diagrams, the attaching circle of each $2$-handle is represented by the framed knot defining the knot trace $X_m(K)$.
The key point is to construct relative trisections of $D^4$ in such a way that these framed knots can be canceled.
We give examples of trisected $2$-handles below.

\begin{figure}[!tbp]
\centering
\includegraphics[scale=0.23]{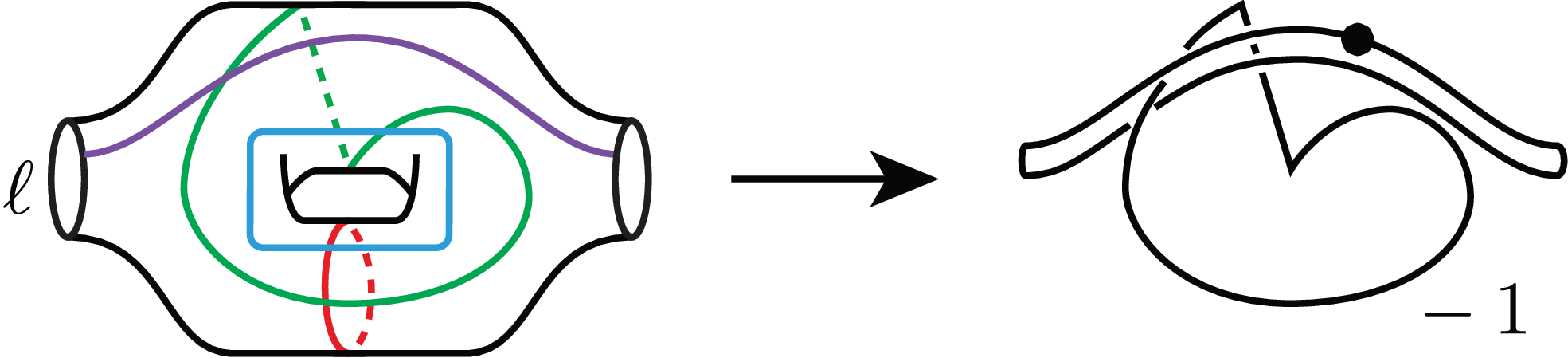}
\caption{$(1,1;0,2)$-relative trisection diagram $\mathcal{D}^{\circ}(\overline{\CC{P^2}}-\Int D^4)$ and the induced handlebody diagram.}
\label{fig:td-cp2bar}
\end{figure}

\begin{example}\label{exm:Ycp2bar}
Let $\mathcal{D}^{\circ}(\overline{\CC{P^2}}-\Int D^4)$ be the diagram shown in the first picture of Figure~\ref{fig:td-cp2bar}.
Capping off the distinguished boundary component $\ell$ with a disk yields a $(1,0;0,1)$-relative trisection diagram of $\overline{\CC{P^2}}-\Int D^4$.
Thus $\mathcal{D}^{\circ}(\overline{\CC{P^2}}-\Int D^4)$ is a $(1,1;0,2)$-relative trisection diagram of $D^4$, and the induced handlebody diagram describes the cancellation of an $(-1)$-framed unknot.
%
By Lemma~\ref{lem:glue-Akb-Lambda}, $\mathcal{D}_A\cup_{\ell_i\sim\ell}\mathcal{D}^{\circ}(\overline{\CC{P^2}}-\Int D^4)$ represent $Y_i(U,-1)$.
\end{example}

\begin{figure}[!tbp]
\centering
\includegraphics[scale=0.23]{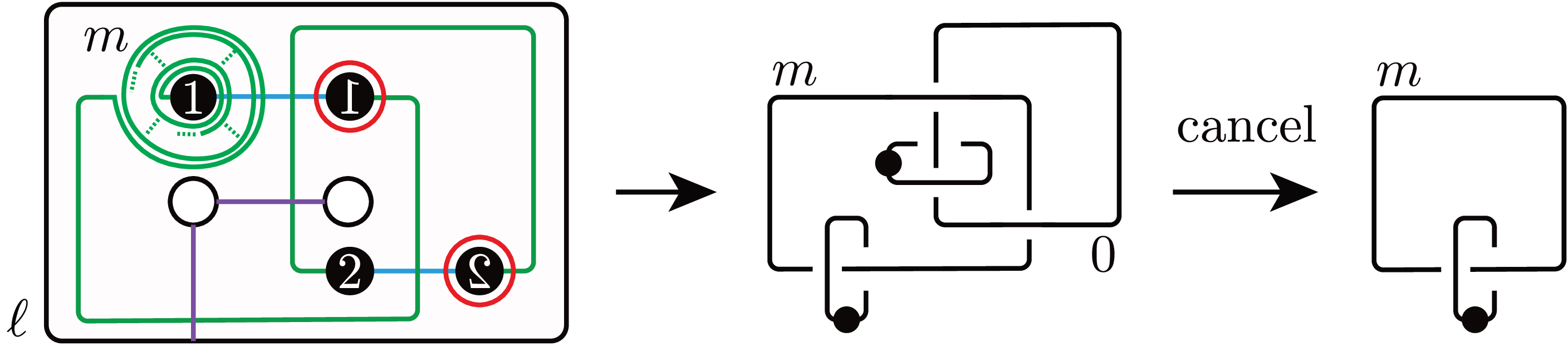}
\caption{$(2,2;0,3)$-relative trisection diagram $\mathcal{D}^{\circ}(U,m)$ and the induced handlebody diagrams.}
\label{fig:td-X_m(U)}
\end{figure}

\begin{example}\label{exm:Yunknot}
Let $\mathcal{D}^{\circ}(U,m)$ be the diagram shown in the first picture of Figure~\ref{fig:td-X_m(U)}.
Capping off the distinguished boundary component $\ell$, which is the outer boundary component in the figure, with a disk yields a $(2,1;0,2)$-relative trisection diagram of $X_m(U)$; see Proposition~4.1 of~\cite{Tak26x}.
Thus $\mathcal{D}^{\circ}(U,m)$ is a $(2,2;0,3)$-relative trisection diagram of $D^4$, and the induced handlebody diagram describes the cancellation of an $m$-framed unknot; see the subsequent diagrams in Figure~\ref{fig:td-X_m(U)}.
By Lemma~\ref{lem:glue-Akb-Lambda}, $\mathcal{D}_A\cup_{\ell_i\sim\ell}\mathcal{D}^{\circ}(U,m)$ represent $Y_i(U,m)$.
The relative trisection diagrams shown in Figures~\ref{fig:td-exoticY1} and~\ref{fig:td-exoticY2} represent $\mathcal{D}_A \cup_{\ell_1\sim\ell} \mathcal{D}^{\circ}(U,-m)$ and $\mathcal{D}_A\cup_{\ell_2\sim\ell}\mathcal{D}^{\circ}(U,-m)$, respectively.
\end{example}

\begin{figure}[!tbp]
\centering
\includegraphics[scale=0.23]{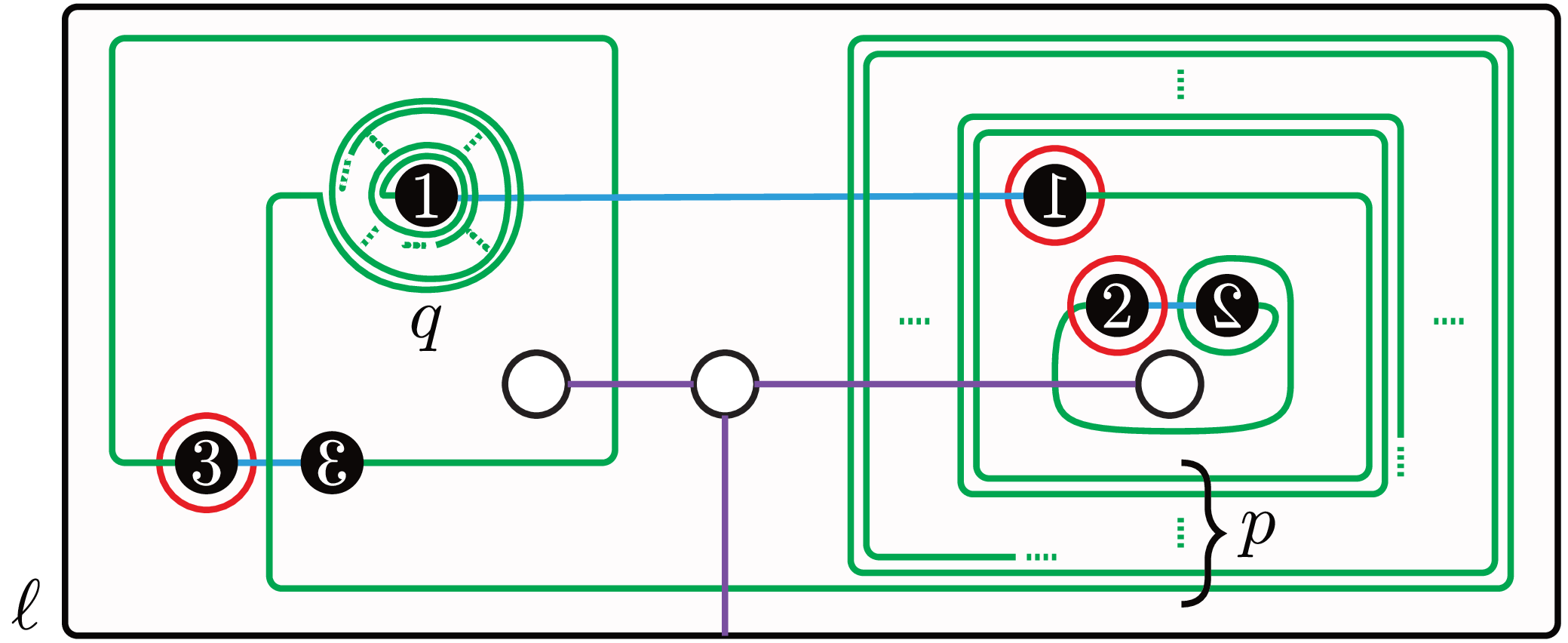}
\caption{$(3,3;0,4)$-relative trisection diagram $\calD^{\circ}(T_{p,p-1},m)$, where $q=m-p^2-p$ is the number of right-handed Dehn twists.}
\label{fig:td-torus}
\end{figure}

\begin{figure}[!tbp]
\centering
\includegraphics[scale=0.23]{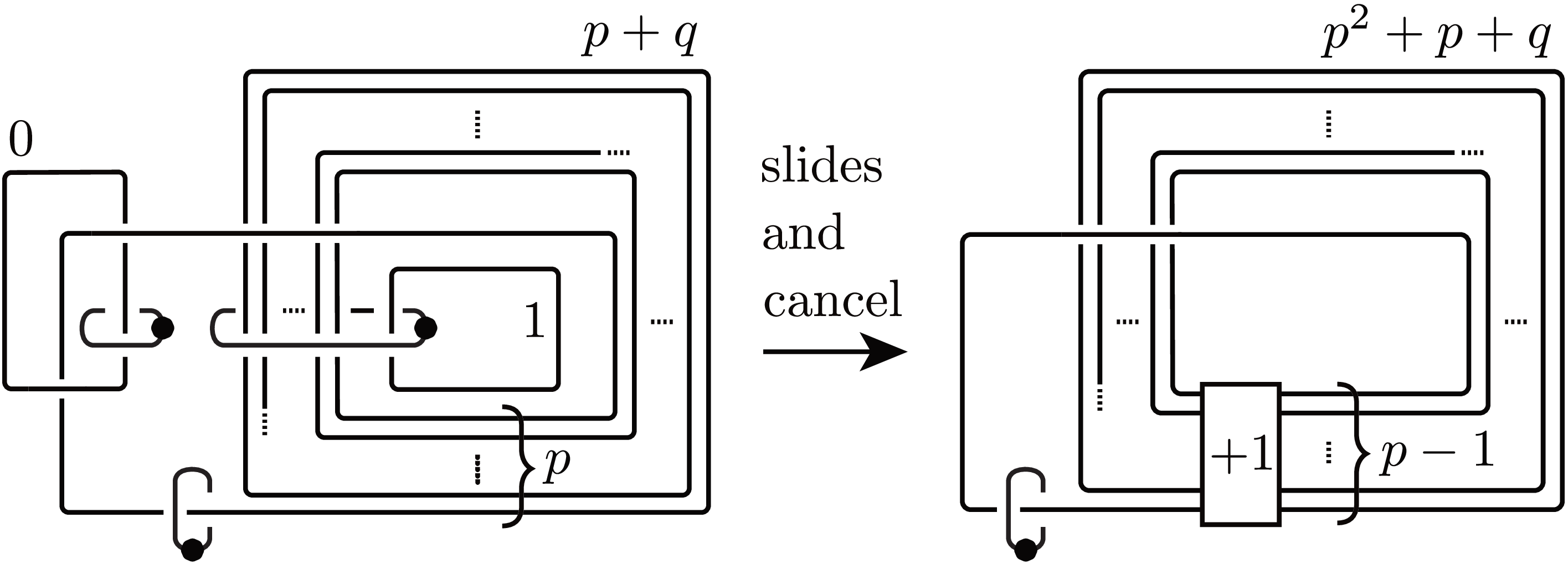}
\caption{The handlebody diagram induced by $\calD^{\circ}(T_{p,p+1},m)$, where $q=m-p^2-p$.}
\label{fig:kd-torus}
\end{figure}

\begin{example}\label{exm:Ytorus}
Let $\mathcal{D}^{\circ}(T_{p,p+1},m)$ be the diagram shown in Figure~\ref{fig:td-torus}.
Capping off the distinguished boundary component $\ell$, which is the outer boundary component in the figure, with a disk yields a $(3,2;0,3)$-relative trisection diagram of $X_m(T_{p,p+1})$; see Proposition~4.2 of \cite{Tak26x}.
Thus $\mathcal{D}^{\circ}(T_{p,p+1},m)$ is a $(3,3;0,4)$-relative trisection diagram of $D^4$, and the induced handlebody diagram describes the cancellation of an $m$-framed torus knot $T_{p,p+1}$; see Figure~\ref{fig:kd-torus}.
By Lemma~\ref{lem:glue-Akb-Lambda}, $\mathcal{D}_A\cup_{\ell_i\sim\ell}\mathcal{D}^\circ(T_{p,p+1},m)$ represent $Y_i(T_{p,p+1},m)$.
\end{example}

\subsection{Proof of the main theorem}

\begin{proof}[Proof of Theorem~\ref{main:exotictris}]
By Examples~\ref{exm:Yunknot} and~\ref{exm:Ytorus}, for any integers $m$ and $p$, the 4-manifolds $Y_1(T_{p,p-1},m)$ and $Y_2(T_{p,p-1},m)$ admit splittable relative trisection diagrams satisfying the conditions of the theorem.
By Theorem~\ref{thm:exoticAY}, these $4$-manifolds are exotic if $p=2$ and $m\leq -1$, or if $p\geq 3$ and $m\leq p^2-3p+1$.
\end{proof}

\begin{proof}[Proof of Cororally~\ref{cor:obd-exotic}]
By Theorem~\ref{thm:exoticAY}, $Y_1(U,-1)$ and $Y_2(U,-1)$ are exotic.
By Example~\ref{exm:Ycp2bar}, each $Y_i(U,-1)$ admits a splittable relative trisection diagram $\mathcal{D}_A\cup_{\ell_i\sim\ell} \mathcal{D}^{\circ}(\overline{\CC{P^2}}-\Int D^4)$.
The diagram $\mathcal{D}^{\circ}(\overline{\CC{P^2}}-\Int D^4)$ induces the open book with annular pages and monodromy given by the right-handed Dehn twist along the core circle.
By Lemma~\ref{lem:srtd-obd}, the open book induced by $\mathcal{D}_A\cup_{\ell_i\sim\ell} \mathcal{D}^{\circ}(\overline{\CC{P^2}}-\Int D^4)$ is obtained by gluing the open books induced by the two constituent pieces.
Consequently, the monodromy of the open book on $\partial Y_i(U,-1)$ is obtained from the monodromy induced by $\mathcal{D}_A$ by composing with the right-handed Dehn twist along $\ell_i$.
\end{proof}

\begin{remark}
For each integer $m\leq -5$, relative trisections for the exotic pair $Y_1(T_{2,3},m)$ and $Y_2(T_{2,3},m)$ can also be obtained from certain Lefschetz fibrations.
Ukida~\cite{Uki16} constructed a positive allowable Lefschetz fibration on each $Y_i(T_{2,3},m)$, with regular fiber $\Sigma_{0,5-m}$ and $5-m$ vanishing cycles.
Since an achiral Lefschetz fibration over $D^2$ can be converted into a relative trisection~\cite{CasGayPin18a}, 
Ukida's PALFs yield relative trisection diagrams of $Y_i(T_{2,3},m)$ with parameters $(5-m,4-m;0,5-m)$.
However, the resulting diagrams are not splittable.
\end{remark}

\begin{remark}
Genus-3 relative trisection diagrams of the Akbulut cork have been constructed in \cite{Tak25a} and \cite{Tak25b}.
However, these diagrams do not give rise to splittable trisections of exotic pairs.
In other words, gluing a relative trisection diagram along boundary components of these diagrams does not produce a pair of homeomorphic $4$-manifolds.
\end{remark}


\section*{Acknowledgements}
The author was partially supported by JSPS KAKENHI Grant Number 26K24497.

\section*{AI disclosure}
The author used ChatGPT (OpenAI) to improve the English writing of this paper.
All mathematical statements, proofs, and figures were prepared by the author, who takes full responsibility for the content of this paper.



\end{document}